\documentclass{article} 
\usepackage{amssymb}
\usepackage{amsmath}
\usepackage{caption2}
\input epsf.sty

\begin{document}

\title{Solution of the Hurwitz problem for Laurent polynomials
}  
\author{F. Pakovich}
\date{}

\maketitle


\def\be{\begin{equation}}
\def\ee{\end{equation}}
\def\bs{$\square$ \vskip 0.2cm}
\def\d{{\rm d}} 
\def\D{{\rm D}} 
\def\I{{\rm I}} 
\def\C{{\mathbb C}} 
\def\N{{\mathbb N}} 
\def\P{{\mathbb P}}
\def\Z{{\mathbb Z}}
\def\R{{\mathbb R}} 
\def\ord{{\rm ord}}

\def\e{\eqref}
\def\phi{{\varphi}}
\def\v{{\varepsilon}} 
\def\O{{\Omega}}
\def\deg{{\rm deg\,}} 
\def\div{{\rm div\,}} 
\def\Det{{\rm Det}}
\def\dim{{\rm dim\,}} 
\def\Ker{{\rm Ker\,}} 
\def\Gal{{\rm Gal\,}}
\def\St{{\rm St\,}} 
\def\exp{{\rm exp\,}} 
\def\cos{{\rm cos\,}} 
\def\diag{{\rm diag\,}} 
\def\GCD{{\rm GCD }}
\def\LCM{{\rm LCM }}
\def\mod{{\rm mod\ }}
\def\S{\C\P^{1}}
\def\O{\Omega}
\def\sk{{\rm sk}}
\def\c{{\rm c}}

\def\bp{\begin{proposition}}
\def\ep{\end{proposition}}
\def\bt{\begin{theorem}}
\def\et{\end{theorem}}
\def\be{\begin{equation}}
\def\bee{\begin{equation*}}
\def\l{\label}
\def\la{\lambda}
\def\m{\mu}
\def\ee{\end{equation}}
\def\eee{\end{equation*}}
\def\bl{\begin{lemma}}
\def\el{\end{lemma}}
\def\bc{\begin{corollary}}
\def\ec{\end{corollary}}
\def\pr{\noindent{\it Proof. }}
\def\note{\noindent{\bf Note. }}
\def\bd{\begin{definition}}
\def\ed{\end{definition}}

\newtheorem{theorem}{Theorem}[section]
\newtheorem{lemma}{Lemma}[section]
\newtheorem{definition}{Definition}[section]
\newtheorem{corollary}{Corollary}[section]
\newtheorem{proposition}{Proposition}[section]

\section{Introduction}

Let $f\,:\, S^2 \rightarrow S^2$ be an $n$-fold branched covering or equivalently a rational function on the Riemann sphere and $z_1, z_2, ... ,z_q \in S^2$
be its branching points (i.e. points $z\in S^2$ for which $f^{-1}\{z\}$ contains less than $n$ points).
Then for each $i,$ $1\leq i \leq q,$ the set
$\Pi_i=\{a_{i,1}, a_{i,2},\, ... \, ,a_{i,p_i}\}$ of local degrees of $f$ at points of $f^{-1}\{z_i\}$ is a 
partition of $n$. Furthermore, it follows from the Riemann-Hurwitz formula that 
\be \l{rh}  
\sum_{i=1}^q p_i=(q-2)n+2.
\ee The collection $\Pi=\{\Pi_1,\, ... \, ,\Pi_q\}$ is 
called the branch datum of $f.$  
In this paper we investigate the following existence problem for rational 
functions:
for a given collection $\Pi$ of partitions $\Pi_i=\{a_{i,1}, a_{i,2},\, ... \, ,a_{i,p_i}\},$    
$1\leq i \leq q,$
of a number $n$ such that \eqref{rh} holds
to define whether there exists a rational function $f$ 
for which $\Pi$ is the branch datum.

The existence
problem for rational 
functions is a particular case of the existence problem for branched coverings
$f\,:\, N \rightarrow M$
between closed Riemann surfaces which goes back to Hurwitz \cite{hu}. 
This problem was studied by many authors (see e.g. \cite{b}-\cite{h}, \cite{pp}, \cite{pp1})
and essentially remains open only for the case when $M=S^2$.
Namely, the results obtained in \cite{eks}, \cite{e}, \cite{h} imply that 
if $\chi(M)\leq 0$ then natural necessary conditions, including 
the Euler characteristic and the orientability of $M$ and $N$ as well as the degree of $f$  
and its local degrees at the branching points, are also sufficient. On the other hand,
if $N$ is the projective plane then Proposition 2.3 of \cite{eks} reduces the problem to the case 
$N=S^2.$ 

In contrast to the case $\chi(M)\leq 0$ if $M=S^2$ then necessary conditions above (which reduce in this case to 
the Riemann-Hurwitz formula) in general are know to be not sufficient.
For example, the collection 
$\{2,2\},$ $\{2,2\},$ $\{3,1\}$ is compatible with \eqref{rh} nevertheless it can not be the branch datum of a
rational function (see \cite{eks}, Corollary 6.4 and the Theorem below). A survey of known results and techniques related to the 
existence problem for branched coverings can be found in \cite{pp}.

The existence problem for branched coverings is closely related to 
the problem of
enumeration of equivalence classes of covering with prescribed branch datum
posed by Hurwitz \cite{hu}. Note that this last problem 
in a sense can be solved using the representation theory of the symmetric group (see \cite{m1}, \cite{m2}),
nevertheless 
the corresponding formulas are usually too complicated to be calculated exactly. 
In particular, an explicit criterion which permits to 
define whether a collection of partitions 
is the branch datum for at least one rational function 
does not exist.

An important particular case when the answer to the existence problem for rational 
functions is known is the one when 
the collection $\Pi$ contains a partition
consisting of a single element. It was shown in 
\cite{t} (see also \cite{eks}, \cite{kz}, \cite{lz}) that for any such a
collection necessary condition \eqref{rh} is also sufficient for the existence
of a rational function for which $\Pi$ is 
the branch datum. Note that the requirement imposed on $\Pi$ implies that this rational function is
equivalent to a polynomial.

Since in view of the remark made above 
the polynomial case seems to be rather special 
the following particular case of the existence problem for rational 
functions, in a sense the simplest possible after the polynomial one, 
is of interest: to describe the collections of partitions, containing a partition 
$\Delta$ consisting of {\it two} elements, 
which are branch date of rational functions. Clearly, this problem 
is essentially equivalent to the existence problem 
for Laurent polynomials.
To our knowledge the only results relevant 
to this problem are: Proposition 5.3 of \cite{eks} which provides the solution of the general existence problem for coverings
in the case when $\Delta=\{1,n-1\}$, Theorem 1.1
of \cite{pp1} which solves the existence problem for Laurent polynomials in the case when $\Delta=\{2,n-2\}$
under the additional assumption that $q=3$, and 
Corollary 6.4 of \cite{eks} which states that a Laurent polynomial with ramification  
$\{2,2,\, ...\, ,  2\},$\ $\{2,2,\, ...\, ,  2\},$\ 
$\{s,n-s\}$ exists if and only if $s=n/2.$

In this paper we provide the complete solution of the existence problem for Laurent polynomials. To formulate our result explicitly let us 
introduce the following notation.
Say that a collection $\Pi$ of $q$ partitions $\Pi_i
=\{a_{i1},a_{i2}, ... , a_{ip_i}\},$ 
$1\leq i \leq q,$ of a number $n$ is an $(n,q)$-passport if the numbers $p_i,$ $1\leq i \leq q,$ are less
than $n$ and satisfy \eqref{rh}. Say that a passport $\Pi$ is realizable if $\Pi$ 
coincides with the branch datum of a rational function. Finally, say that a passport $\Pi$ is a Laurent passport if $p_q=2.$   
Under this notation our main result is the following theorem. 
\vskip 0.2cm
\noindent{\bf Theorem.} {\it
Any Laurent passport $\Pi$ for which $q>3$ is realizable. 
A Laurent passport $\Pi$ for which $q=3$ is realizable if 
and only if $\Pi$ is distinct from the triplets listed below:
\vskip 0.1cm
\noindent 1)\ 
$\{l,l,\, ...\,  ,l\}$,\ $\{1,1,\, ...\, ,1, d\},$\ 
$\{s,n-s\},$ where $d\geq 3,$ $l\geq 2,$ $s\geq 1,$  
$s\equiv 0\ \mod l,$ 
\vskip 0.15cm
\noindent 2)\
$\{2,2,\, ...\, ,  2\},$\ $\{2,2,\, ...\, ,  2\},$\ 
$\{s,n-s\},$ where $s\geq 1,$ $s\neq n/2,$ 
\vskip 0.15cm
\noindent 3)\
$\{2,2, \, ... \, , 2\},$\ $\{1,1,\, ... \, , 1,d-1,d\},$\ 
$\{2d-3,n-2d+3\},$ where $d\geq 3,$ \vskip 0.2cm
\noindent 4)\
$\{2,2, \, ... \, , 2\},$\ $\Pi_2=\{1,1,\, ... \, , 1,d,d\},$\ 
$\Pi_3=\{2d-3,n-2d+3\},$ 
where $d\geq 3,$
\vskip 0.15cm
\noindent 5)\
$\{2,2, \, ... \, , 2\},$ $\{1,1,\, ... \, , 1,d,d\},$ 
$\{2d-1,n-2d+1\},$ where $d\geq 3,$
\vskip 0.15cm
\noindent 6)\
$\{2,2, \, ... \, , 2\},$\ $\Pi_2=\{1,2,2,\, ... \, ,2 
,3\},$\  $\Pi_3=\{n/2,n/2\},$ 
\vskip 0.15cm
\noindent 7)\
$\{2,2,2,2,2,2\},$ $\{1,1,1,3,3,3\},$ 
$\{6,6\}.$
}

\pagebreak

Our approach to the existence problem for rational functions is based on 
a one-to-one correspondence between equivalence classes of $n$-fold branched coverings $f\,:\, S^2 \rightarrow S^2$ with 
branching points $c_1, c_2,\, ...\, ,c_q$ and  
equivalence classes of so called {\it planar $(n,q)$-constellations}
(see \cite{lz} and section 2 below).
Roughly speaking
a planar $(n,q)$-constellation is a connected planar graph $\Gamma$ obtained by gluing together $n$ copies
of a planar $(q-1)$-gone with numerated vertices along vertices 
with equal numbers, and
to a covering with the branch date $\Pi=\{\Pi_1,\, ... \, ,\Pi_q\}$ corresponds a constellation
for which $\Pi_i,$ $1\leq i \leq q-1,$ coincides with the collection of valencies of 
vertices of $\Gamma$ with the number $i$ while $\Pi_q$ is related with  
the collection of valencies of faces of $\Gamma$.
The correspondence between coverings and constellations reduces the  
existence problem for rational functions 
with prescribed branch data   
to the existence problem
for constellations
with prescribed valency data
and in this paper we will consider the existence problem in this
purely combinatorial setting.

Note that in the case when $q=3$ constellations are simply bicolored planar graphs
that is the graphs whose vertices can be colored by two colors so that adjacent vertices have different colors.
Such graphs, also called ``dessins d'enfants", are 
closely related to Galois theory (see e.g. \cite{lz} and the bibliography there) and for this reason appear in a large number of recent papers.
In general case however constellations have more subtle combinatorial structure
and  
one of the objectives of this paper is to develop some combinatorial techniques
to work with constellations  
in order to make these beautiful combinatorial objects  
useful for the questions like the existence problem.
Note also that since with appropriate modifications the correspondence 
above extends to a correspondence between coverings $f\,:\, N \rightarrow S^2$, where
$N$ is any closed Riemann surface (which is necessarily orientable) and the 
corresponding graphs embedded in $N$, our method in principle is applicable
for such coverings too.

The paper is organized as follows. In the second section we recall the 
correspondence between constellations and coverings and introduce the notation.
Besides, we prove two lemmas which we will often use in the following. 
In the third section we develop the necessary techniques and give
the constructive proof of the main theorem in the case $q>3.$
Finally, in the fourth section we separately analyse the case $q=3$ 
which turns out to be essentially different from the general one.

\section{Preliminaries and notation}
\subsection{Constellations and coverings}
In this subsection we recall the 
correspondence between constellations and coverings.
For more information and other versions of the definition of a constellation we refer the reader to \cite{lz}.  

{\it A $q$-star} is a connected planar graph $S$, consisting of one vertex of valency $q,$ $q$ vertices of valency $1,$ and $q$ edges, 
such that the vertices of valency $1$ are numerated in the counterclockwise direction 
with respect to the natural cyclic ordering induced by the embedding of $S$ (see Fig. \ref{c1},a).
{\it A planar $(n,q)$-constellation} $\Gamma$ is a connected planar graph obtained by gluing together 
$n$ copies of a $q-1$-star along their numerated vertices 
with equal numbers
(see
Fig. \ref{c1},b). 
\begin{figure}[h]
\renewcommand{\captionlabeldelim}{.}
\medskip
\epsfxsize=11truecm
\centerline{\epsffile{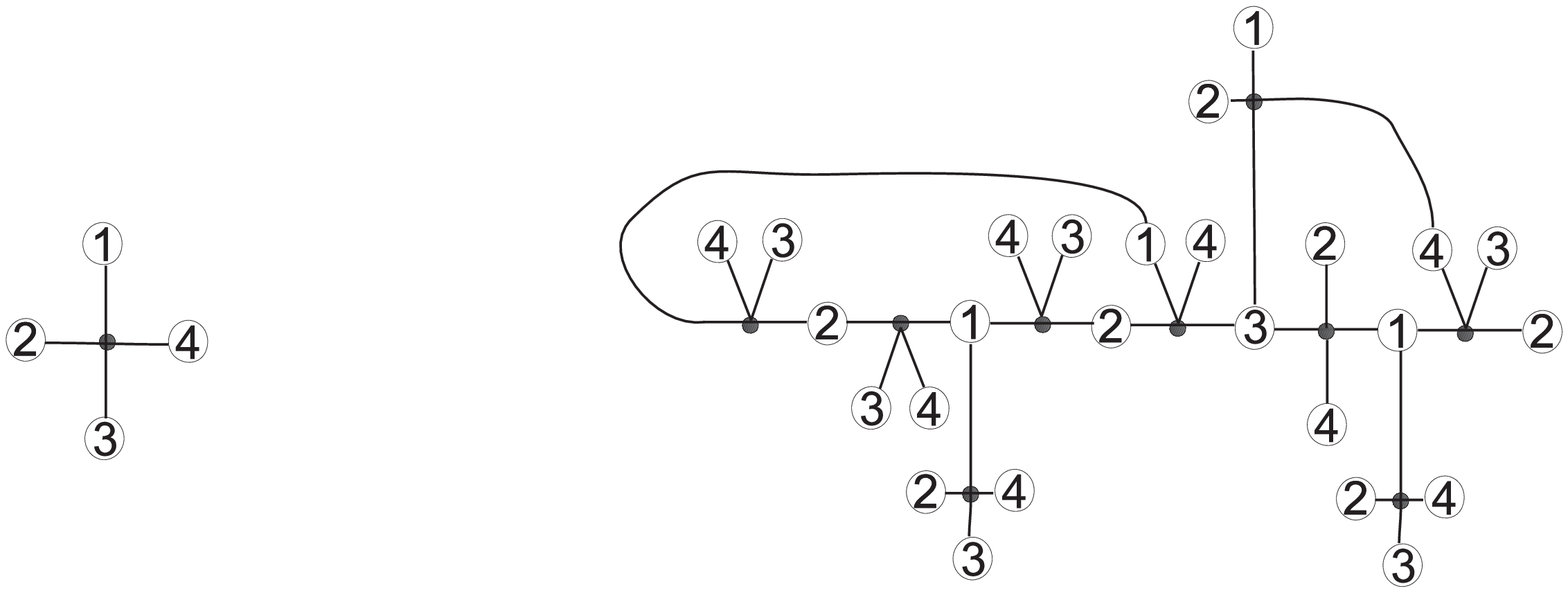}}
\caption{}\l{c1}
\medskip
\end{figure}We will suppose additionally that for each $i,$ $1\leq i \leq q-1,$ 
the graph $\Gamma$ contains a vertex with number $i$ whose valency is $\geq 2$
and that the number of poles of $\Gamma$ is less that $n$.   
Two planar constellations $\tilde\Gamma$ and $\Gamma$ are called {\it equivalent} if 
$\tilde\Gamma=h(\Gamma)$,
where $h\,:\, S^2 \rightarrow S^2$ is an orientation preserving homeomorphism which
preserves the numbers of vertices. Since in this paper we will work only with planar constellations in the following we will omit the word ``planar''.
Note that if we traverse a face of a constellation $\Gamma$ then the numbers of numerated vertices appear in the cyclic or inverse cyclic order and between any two consecutive numerated vertices there is 
exactly one non-numerated vertex. In particular, the valency of each face  
of $\Gamma$ is divisible by $2(q-1).$ 

The numerated vertices of a constellation $\Gamma$ with number $i,$ $1\leq i \leq q-1,$ are called 
{\it $i$-vertices} of $\Gamma$ and the collection of valencies of $i$-vertices of $\Gamma$ is denoted by $\Gamma_i=\{a_{i,1},a_{i,2}, ... , a_{i,p_i}\}$. 
By $\Gamma_{q}=\{a_{q,1},a_{q,2}, ... , a_{q,p_q}\}$ 
we will denote the collection of valencies of faces of $\Gamma$ divided by $2(q-1).$
Note that in view of the remark above for any $i,$ $1\leq i \leq q-1,$
the number $a_{q,j},$ $1\leq j \leq p_q,$ equals the number of appearances of $i$-vertices when traversing the corresponding face.
We will call the collection $\Gamma_1, \Gamma_2,\, ... \, , \Gamma_q$ the 
{\it valency datum} of the constellation $\Gamma.$ For example, for a $(9,5)$-constellation shown on  
Fig. \ref{c1},b its valency datum is $\Gamma_1=\{1,2,3,3\},$ $\Gamma_2=\{1,1,1,1,1,2,2\},$ 
$\Gamma_3=\{1,1,1,1,1,1,3\},$ $\Gamma_4=\{1,1,1,1,1,1,1,2\},$ 
$\Gamma_5=\{1,2,6\}.$ 

Since each star of a constellation $\Gamma$ is adjacent to a unique $i$-vertex of $\Gamma$  
each collection $\Gamma_i=\{a_{i,1},a_{i,2}, ... , a_{i,p_i}\},$ $1\leq i \leq q-1,$ is a partition of $n.$ Furthermore, since the sum of valencies of faces of $\Gamma$ coincides with 
the doubled number of edges of $\Gamma$ the collection $\Gamma_q=\{a_{q,1},a_{q,2}, ... , a_{q,p_q}\}$ also is a partition of $n.$ 
Notice that the additional requirement made in the definition of a constellation is
equivalent to the requirement that the numbers $p_i,$ $1\leq i \leq q,$ are less than $n$. 
Finally, observe that Euler's formula implies that the numbers $p_i,$ $1\leq i \leq q,$ satisfy \eqref{rh}. 

Starting from an $n$-fold branched covering $f\,:\, S^2 \rightarrow S^2$ with $q$ 
branching points $c_1, c_2,\, ...\, ,c_q$ and the branch datum $\Pi=\{\Pi_1,\, ... \, ,\Pi_q\}$
we can obtain an $(n,q)$-constellation $\Gamma=\Gamma(f)$ for which $\Gamma_i=\Pi_i,$
$1\leq i \leq q,$ as follows.  
Let $c$ be
a non-branching value of $f(z)$ and $S\subset S^2$ be 
a $q-1$-star joining $c$ with $c_1, c_2,\, ...\, ,c_{q-1}$ 
such that $c_q\in S^2\setminus S$.
Define $\Gamma$ as the preimage of $S$ under the map $f\,:\, S^2 \rightarrow S^2.$
More precisely, define edges of $\Gamma$ as preimages of edges of $S,$  
$i$-vertices of $\Gamma$ as preimages of $c_i,$ $1\leq i \leq q-1,$ and non-numerated vertices
of $\Gamma$ as preimages of $c$
(see
Fig. \ref{c2}).
It is not hard to verify that $\Gamma$ is indeed a constellation and that $\Gamma_i=\Pi_i,$
$1\leq i \leq q,$.  
\begin{figure}[h]
\renewcommand{\captionlabeldelim}{.}
\medskip
\epsfxsize=11truecm
\centerline{\epsffile{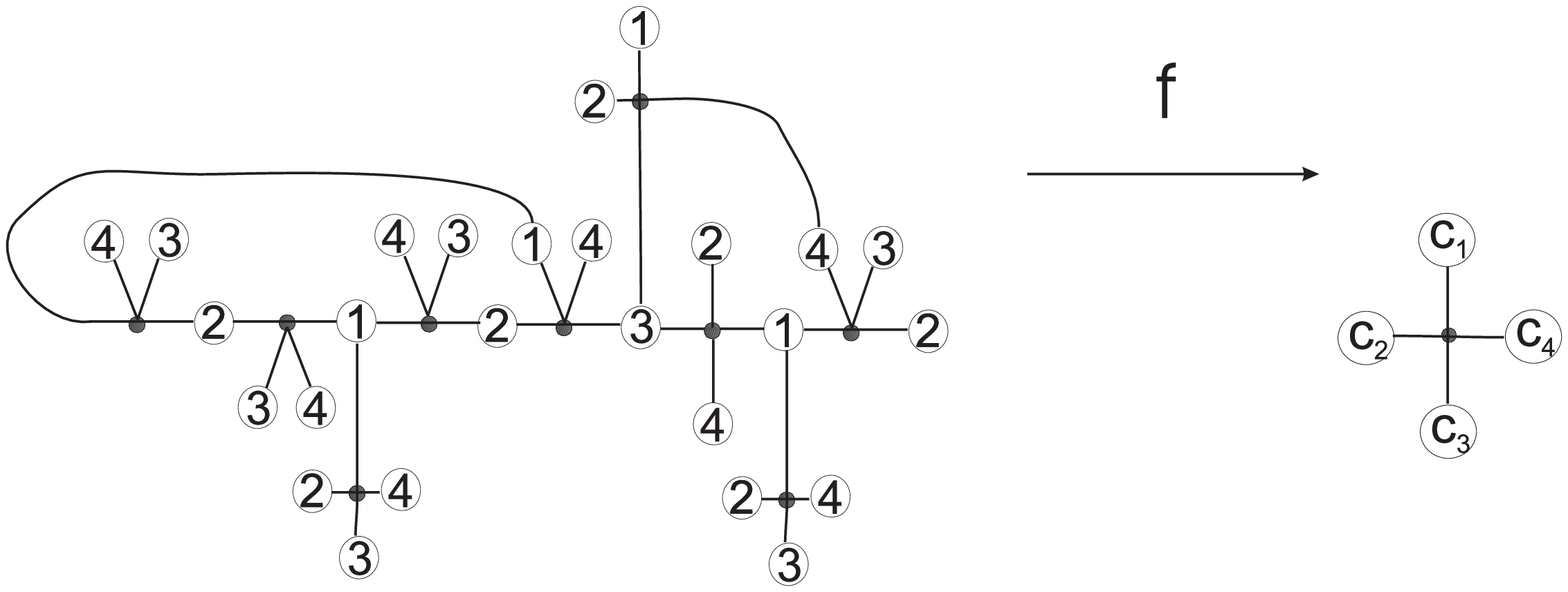}}
\caption{}\l{c2}
\medskip
\end{figure}  

Conversely, if $\Gamma$ is an $(n,q)$-constellation with the valency datum 
$\Gamma_1, \Gamma_2,\, ... \, ,$ $\Gamma_q$ then for any $c_1, c_2,\, ...\, ,c_{q}\in S^2$ 
there exists an $n$-fold branched covering 
$f\,:\, S^2 \rightarrow S^2$ with
branching points $c_1, c_2,\, ...\, ,c_{q}$ and the branch datum $\Pi=\{\Pi_1,\, ... \, ,\Pi_q\}$
such that $\Pi_i=\Gamma_i,$ $1\leq i \leq q.$ 
To construct the covering needed first of all modify the constellation 
$\Gamma$ as follows. 
Encircle each 
star $S_l,$ $1\leq l \leq n,$ of $\Gamma$ with a simple closed curve $\gamma_l$ so that the closure of the 
domain $D_l$ bounded by $\gamma_l$ contains $S_l$, and   
$\gamma_l\cap \Gamma$ consists of numerated vertices of $S_l$ only. Then delete all the edges 
and non-numerated vertices of $\Gamma$ (see
Fig. \ref{c3},a, where this operation is applied to the constellation shown on Fig. \ref{c2}). Clearly, the obtained graph $\Omega$ has a natural two-colored structure 
on his faces. We will color the faces $D_l,$ $1\leq l \leq n,$ by the black color and the rest faces $L_j,$ $1\leq j \leq p_q,$ by the white one. 
\begin{figure}[t]
\renewcommand{\captionlabeldelim}{.}
\medskip
\epsfxsize=11truecm
\centerline{\epsffile{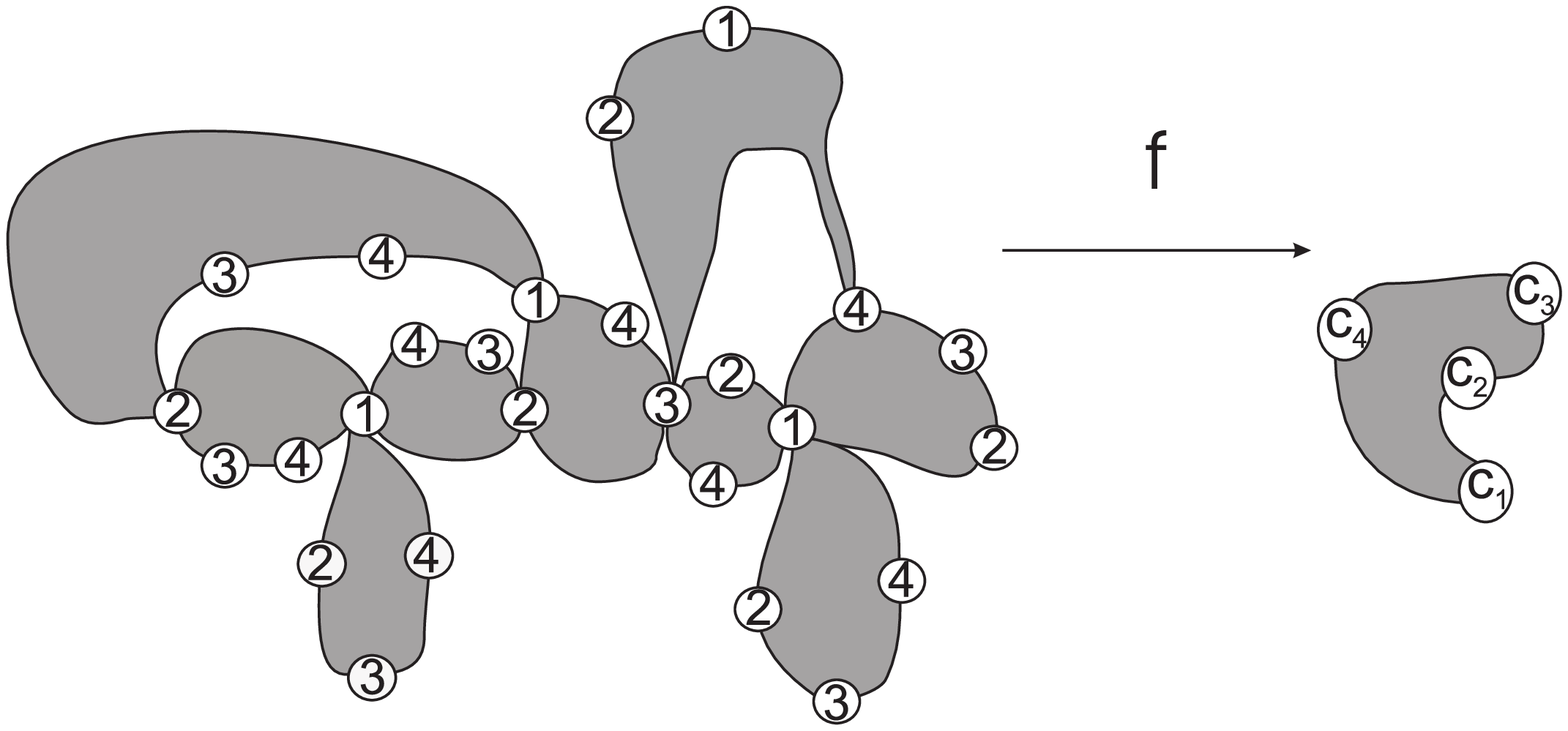}}
\caption{}\l{c3}
\medskip
\end{figure}

Let $\gamma$ be a simple closed curve which passes through $c_1, c_2,\, \ ...\, ,c_{q-1}$ consecutively.
It divides the sphere into two parts. Denote the bounded part by $D$ 
and the unbounded part by $L$ (see Fig. \ref{c3},b, where $D$ (resp. $L$) is colored by black (resp. white) color). Suppose additionally that $\gamma$ is chosen in
such a way that $c_q\in L.$
It is not hard to see that we can define a continuous function 
$f\,:\, S^2 \rightarrow S^2$ which satisfies the following condition: $f$ maps $\bar D_l,$ $1\leq l \leq n,$ on $\bar D$ homeomorphically such that the
$i$-vertex of $\bar D_i$ is mapped on $c_i,$ $1\leq i \leq q,$ while the restriction of $f$ on $L_j,$ 
$1\leq j \leq p_q,$ is a $a_{q,j}$-fold branched covering 
of $L$ with the unique branching point $c_q$ ($f$ on $L_j$ looks like $z^{a_{q,j}}$ on the unit circle). Clearly, 
$f$ is an $n$-branched covering and by construction the valency datum of $\Gamma$ coincides with the branch datum
of $f.$ 

It is easy to check that the correspondence above descends to a one-to-one correspondence
between equivalence classes of $n$-fold branched coverings $f\,:\, S^2 \rightarrow S^2$ with 
branching points $c_1, c_2,\, ...\, ,c_q$ and  
equivalence classes of planar $(n,q)$-constellations. In particular, this implies  
that instead of proving that a covering with a given branch datum exists or does not exist 
it is enough to prove the corresponding fact about constellations.

Notice that $(n,3)$-constellations
are in a one-to-one correspondence with $n$-edged 
bicolored planar graphs. Indeed, it is enough ``to forget'' about non-colored vertices and paint 1-vertices (resp. 2-vertices) by the back (resp. the white) color (see Fig. \ref{c6}). The corresponding rational functions are called Belyi functions and have
very interesting arithmetical properties (see e. g. 
\cite{lz}).

\begin{figure}[h]
\renewcommand{\captionlabeldelim}{.}
\medskip
\epsfxsize=7truecm
\centerline{\epsffile{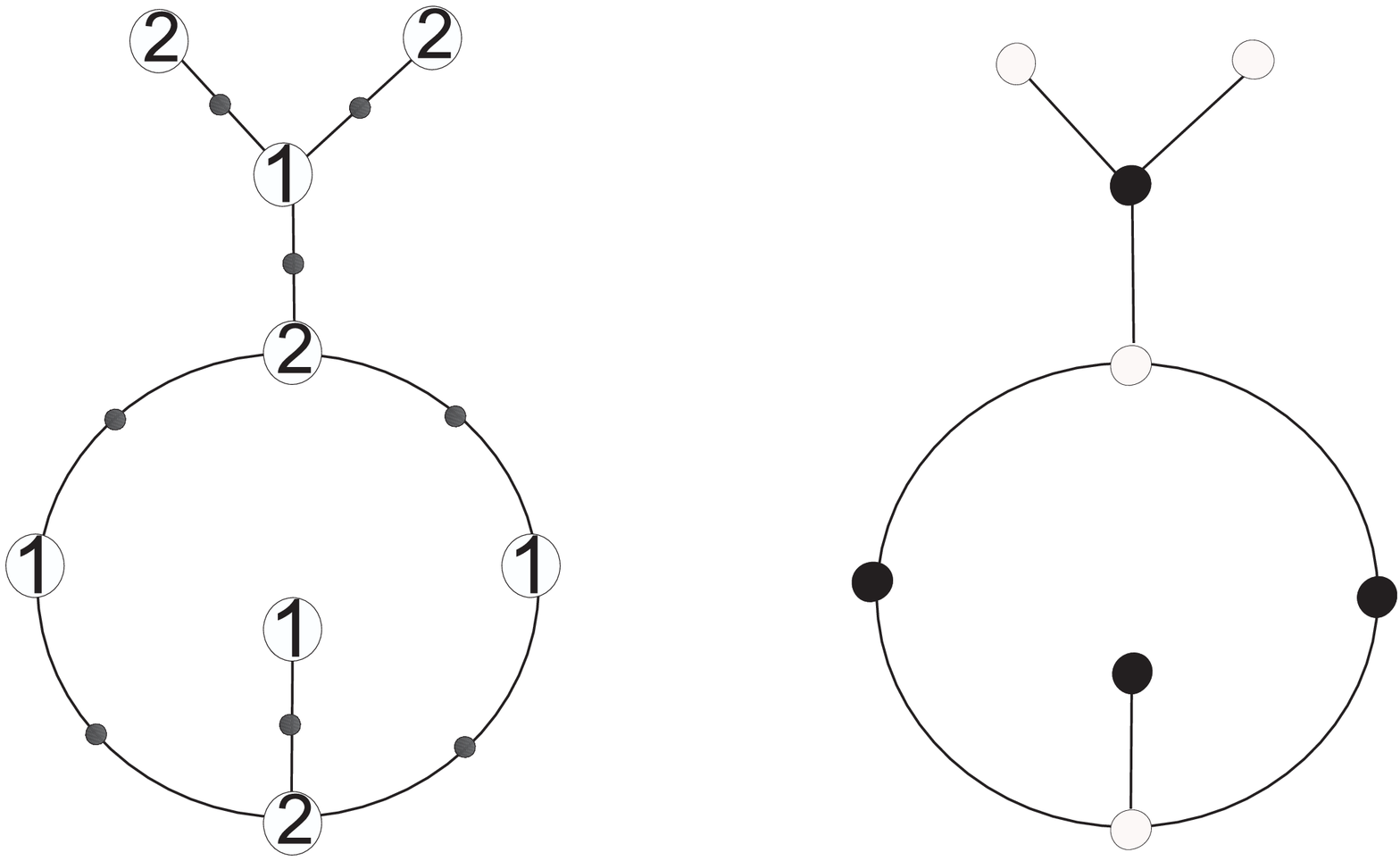}}
\caption{}\l{c6}
\medskip
\end{figure}

\subsection{Constellations with two faces and Laurent passports} 
 
In this subsection we fix notation concerning two-face constellations 
and Laurent passports. Besides, we prove two simple lemmas about constellations and Laurent passports which we will often use in the following.  
\subsubsection{Notation for Laurent passports}

First of all, since for a Laurent $(n,q)$-passport $\Pi$
the partition $\Pi_q=\{s,n-s\}$ essentially depends only on the parameter
$s$ (for given $n$), we will always
indicate only this parameter instead of writing explicitly 
the partition itself. 
Besides, it is convenient 
to denote the number $q-1$ which will appear in most formulas 
by another letter $r$.

Furthermore, for a Laurent passport $\Pi$ we will denote by $q_i$ (resp. $e_i$), $1\leq i \leq r,$ the number of elements of 
$\Pi_i=\{a_{i,1}, a_{i,2},\, ... \, ,a_{i,p_i}\}$ which are greater than 1 (resp. equal 1) and by
$b_{i,1},b_{i,2}, ... , b_{i,q_i},$ $1\leq i \leq r,$ the elements of 
$\Pi_i$ which are greater than 1. Clearly, we
have $e_i+q_i=p_i,$ $1\leq i \leq r,$ and equality \eqref{rh} reduces to the equality  
\be \l{rh1}  
\sum_{i=1}^r p_i=(r-1)n.
\ee
{\it To be definite we will always assume that 
$b_{i,1}\leq b_{i,2}\leq ... \leq b_{i,q_i},$ $1\leq i \leq r,$ and}
$q_1\geq q_2 \geq ... \geq q_{r}.$ 

\subsubsection{Notation for constellations with two faces}

First of all notice that although a constellation is an object embedded in $S^2$ 
all our pictures will be plane. In view of this fact we will use the following notation. For a pictured two-face constellation 
a bounded (resp.
an unbounded) face of $\Gamma$ is called {\it an interior} (resp.  
{\it an exterior}) face of $\Gamma$. To lighten notation the corresponding number $a_{q,i}\in\Gamma_q,$ $i=1,2,$ is denoted by $i(\Gamma)$ (resp. $e(\Gamma)$). 
\begin{figure}[t]
\renewcommand{\captionlabeldelim}{.}
\medskip
\epsfxsize=8truecm
\centerline{\epsffile{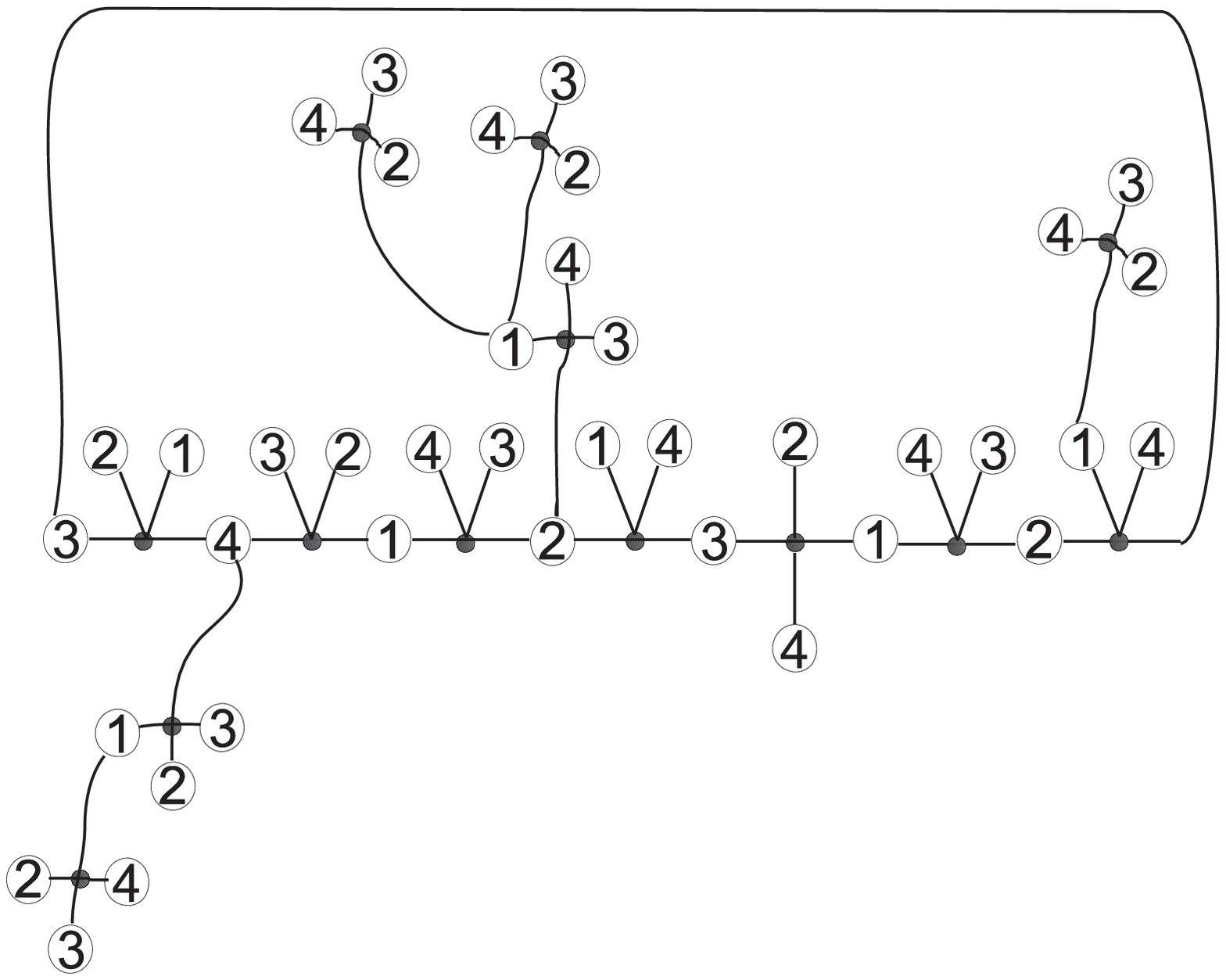}}
\caption{}\l{c7}
\medskip
\end{figure}

Furthermore, a union of all stars of a two-face constellation $\Gamma$
which have an edge adjacent to both faces of $\Gamma$ is called 
{\it a skeleton} of $\Gamma$ and is denoted by $\sk(\Gamma).$ The graph obtained from 
$\sk(\Gamma)$ by removing all vertices of valency $1$, together with adjacent to them edges, and all non-colored vertices is called the cycle of $\Gamma$ and is denoted by $\c(\Gamma).$ 
For example, for the constellation shown on Fig. \ref{c7} the corresponding skeleton and cycle 
are shown on Fig. \ref{ske}.

Let $v$ be a numerated 
vertex of $\Gamma$ adjacent to a star which belongs to $\sk(\Gamma)$.
A subconstellation $\lambda$ of $\Gamma$ such that $\lambda$ contains $v$,
$\lambda\setminus v$ 
belongs to the 
bounded (resp. the unbounded) part of $S^2\setminus \sk(\Gamma)$,  
and $\Gamma\setminus \lambda$ is connected
is called {\it an interior} (resp. {\it an exterior}) {\it branch} of $\Gamma$ growing from $v.$
The number of stars of a branch $\lambda$ is called {\it the weight} of $\lambda$ and is denoted by $\vert\lambda\vert.$ 
For example, the constellation shown on Fig. \ref{c7} has one exterior branch of weight $2$ and two interior branches whose weights are $1$ and $3.$ A constellation $\Gamma$ which 
does not have interior branches is called {\it a sunflower}.

It is convenient to use for two-face constellations the notation 
similar to the one for Laurent passports. So, for a two-face $(n,q)$-constellation $\Gamma$
we will denote by $r$ the number $q-1,$
by $q_i$ (resp. $e_i$) the number of elements of 
$\Gamma_i=\{a_{i,1}, a_{i,2},\, ... \, ,a_{i,p_i}\},$ $1\leq i \leq r,$ which are greater than 1 (resp. equal 1),
and by $b_{i1},b_{i2}, ... , b_{iq_i},$ $1\leq i \leq r,$ the elements of 
$\Gamma_i$ which are greater than 1.
To avoid any confusion in case of necessity
we will write in parenthesis to which passport or constellation these quantities 
and the parameters $n,r$ are related. Clearly, formula \eqref{rh1} holds also 
for two-face $(n,q)$-constellations.

Since in the rest of this paper we will deal 
only with passports which are Laurent and with constellations which are two-faced 
in the following we will omit the corresponding adjectives. 
\begin{figure}[h]
\renewcommand{\captionlabeldelim}{.}
\medskip
\epsfxsize=8truecm
\centerline{\epsffile{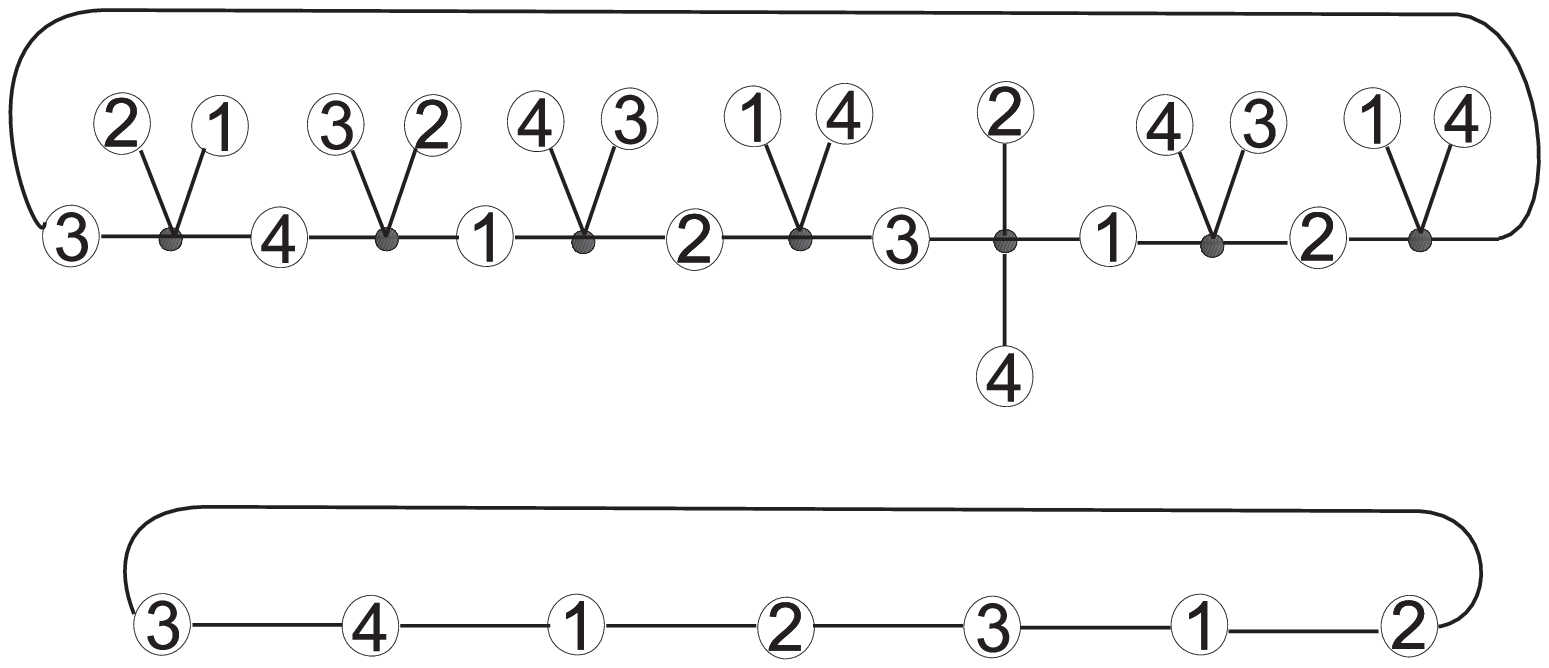}}
\caption{}\l{ske}
\medskip
\end{figure} 

\subsubsection{Two lemmas
}

\bl \l{gop} For any passport $\Pi$ or constellation $\Gamma$ we have:
$$\sum_{i=2}^r \sum_{j=1}^{q_i} (b_{i,j}-2)=e_1+q_1-(q_2 + q_3 + ... + q_r)$$
\el

\pr Indeed, $$\sum_{i=2}^{r}\sum_{j=1}^{q_i} (b_{i,j}-2)=\sum_{i=1}^{r} \sum_{j=1}^{q_i} (b_{i,j}-2)-\sum_{j=1}^{q_1} (b_{1,j}-2)=$$
$$=\sum_{i=1}^{r} \sum_{j=1}^{q_i} (b_{i,j}-2)+2q_1-\sum_{j=1}^{q_1} b_{1,j}.$$ On the other hand, 
$$ \sum_{i=1}^r \sum_{j=1}^{q_i} (b_{i,j}-2)-
\sum_{i=1}^{r}e_i=\sum_{i=1}^{r} \sum_{j=1}^{p_i} (a_{i,j}-2)=$$ $$ =nr -2\sum_{i=1}^{r}p_i=
 nr -2(r-1)n=(2-r)n.$$
Therefore,  
$$\sum_{i=2}^r \sum_{j=1}^{q_i} (b_{i,j}-2)=(2-r)n+\sum_{i=1}^{r}e_i+2q_1-\sum_{j=1}^{q_1} b_{1,j}=$$$$=
(2-r)n+\sum_{i=1}^r(e_i+q_i)-\sum_{i=2}^{r}q_i+q_1-\sum_{j=1}^{q_1} b_{1,j}=$$
$$=
(2-r)n+\sum_{i=1}^rp_i-\sum_{i=2}^{r}q_i+q_1-(n-e_1) =$$
$$
=(2-r)n+(r-1)n-\sum_{i=2}^{r}q_i+q_1-(n-e_1)=$$$$=e_1+q_1-(q_2 + q_3 + ... + q_r).$$
\vskip 0.2cm

\bl \l{ed} Let $\Pi$ be a passport and $\Gamma$ be a constellation such that 
$r(\Gamma)=r(\Pi),$ $q_i(\Gamma)=q_i(\Pi),$ $1\leq i \leq r,$ and 
$b_{i,j}(\Gamma)=b_{i,j}(\Pi),$ $1\leq i \leq r,$ $1\leq j \leq q_i.$
Then $\Gamma_i=\Pi_i,$ $1\leq i \leq r.$
\el

\pr Indeed, it follows from lemma \ref{gop} that $e_1(\Gamma)=e_1(\Pi).$ 
Since $b_{1,j}(\Gamma)=b_{1,j}(\Pi),$ $1\leq j \leq q_1,$ this implies that
$\Gamma_1=\Pi_1.$ Therefore, $n(\Gamma)=n(\Pi).$ But then also
$e_i(\Gamma)=e_i(\Pi),$ $2\leq i \leq r,$ and therefore $\Gamma_i=\Pi_i,$ $1\leq i \leq r.$
\vskip 0.2cm

The lemma \ref{ed} implies that in order to prove that a passport $\Pi$ is realizable it 
is enough to find a constellation $\Gamma$ for which $q_i(\Gamma)=q_i(\Pi),$ $1\leq i \leq r,$ and
$b_{i,j}(\Gamma)=b_{i,j}(\Pi),$ $1\leq i \leq r,$ $1\leq j \leq q_i,$
without checking that $n(\Gamma)=n(\Pi)$ and $e_i(\Gamma)=e_i(\Pi),$ $1\leq i \leq r.$
We will often use this fact without mentioning it explicitly.


\section{Passports with $r>2$.} 

\bp \l{chain2} Let $r> 2$ and $q_1 \geq q_2 \geq q_3 \geq  ... \geq q_r>0$ be
integers such that 
$q_1\leq q_2+ q_3 + ... + q_r$. 
Then for any $s,$ $1\leq s \leq q_2+ q_3 + ... + q_r,$ there exists a sunflower $\Omega$ 
such that all numerated vertices of $\Omega$ have valencies $\leq 2,$ 
$q_i(\Omega)=q_i,$ $1\leq i \leq r,$ and $i(\Omega)=s.$ 

\ep

\pr 
We will prove the proposition in three stages. First we will construct a sunflower $\Delta$ for which $q_1(\Delta)=q_2,$ $q_i(\Delta)=q_i,$ $2\leq i \leq r,$ and $i(\Delta)=q_2.$
Then we will construct a sunflower $\Sigma$ such that $q_i(\Sigma)=q_i,$ 
$1\leq i \leq r,$ and $i(\Sigma)=q_1.$ Finally, we will 
construct the sunflower $\Omega.$ 

To construct the sunflower $\Delta$ first dispose $2q_2+q_3+ ... + q_r$ vertices, $q_2$ of which are 1-vertices and $q_i,$ $2\leq i \leq r,$ 
are $i$-vertices, on the circle as follows: place a 1-vertex as the ``first", a 2-vertex as the 
``second'', and so on till a $r$-vertex (we move in the clockwise direction). 
Then place again a 1-vertex and continue as above skipping however those $i$-vertices, $2\leq i \leq r,$ 
which are already out of stock (see Fig. \ref{f1}, where $q_2=3,$ $q_3=2,$ $q_4=1$).
\begin{figure}[h]
\renewcommand{\captionlabeldelim}{.}
\medskip
\epsfxsize=11truecm
\centerline{\epsffile{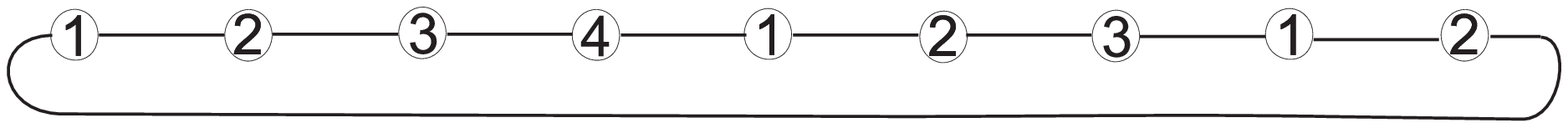}}
\caption{}\l{f1}
\medskip
\end{figure}
Now replace each edge of the obtained graph by a star respecting the vertex numeration as it is
shown on Fig. \ref{f2}. 
\begin{figure}[h]
\renewcommand{\captionlabeldelim}{.}
\medskip
\epsfxsize=11truecm
\centerline{\epsffile{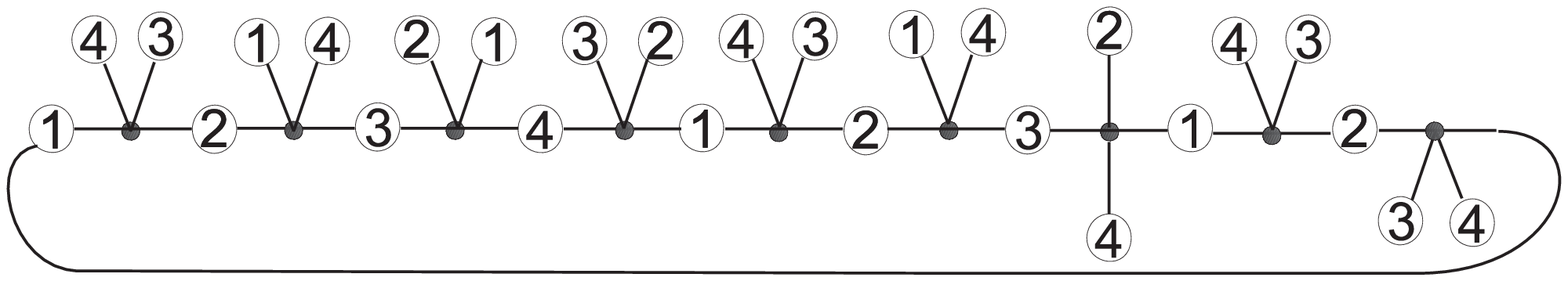}}
\caption{}\l{f2}
\medskip
\end{figure}
Clearly, we obtain a sunflower $\Delta$ for which $q_1(\Delta)=q_2,$ $q_i(\Delta)=q_i,$ $2\leq i \leq r.$ 
Furthermore, the construction implies that 1-vertices of valency 1 can not be adjacent 
to the interior face of $\Delta$. It follows that there are exactly $q_2$ 1-vertices adjacent 
to the interior face of $\Delta$ and hence
the equality $i(\Delta)=q_2$ holds.
 
To construct the sunflower $\Sigma$ modify $\Delta$ 
as follows. Replace any star $S$ of $\Delta$
for which its 1-vertex is of valency 1 (see Fig. \ref{f3},a) 
\begin{figure}[h]
\renewcommand{\captionlabeldelim}{.}
\medskip
\epsfxsize=11truecm
\centerline{\epsffile{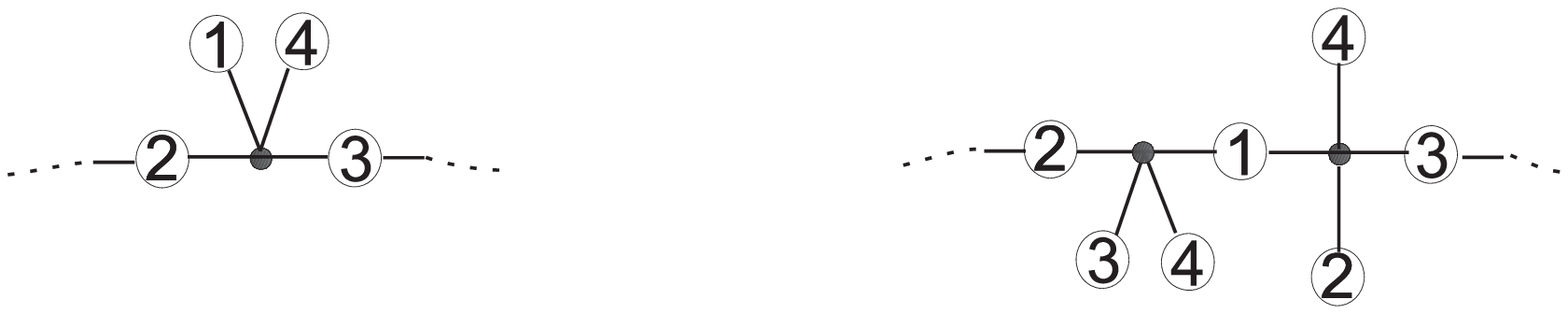}}
\caption{}\l{f3}
\medskip
\end{figure}
by two stars shown on Fig. \ref{f3},b so that to obtain 
a sunflower 
$\tilde \Delta$ such that 
$q_1(\tilde \Delta)=q_1(\Delta)+1$ and $q_i(\tilde \Delta)=q_i(\Delta),$ $2\leq i \leq r$ (see Fig. \ref{f4}, where this operation is applied to the sunflower shown on Fig. \ref{f2}). 
\begin{figure}[h]
\renewcommand{\captionlabeldelim}{.}
\medskip
\epsfxsize=12truecm
\centerline{\epsffile{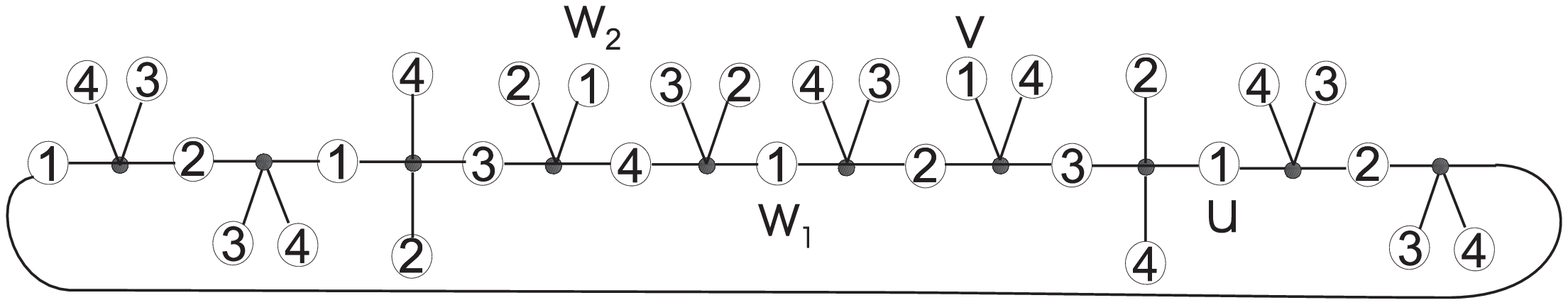}}
\caption{}\l{f4}
\medskip
\end{figure}
Observe that the number of appearances of 
1-vertices when traversing the exterior face of $\tilde \Delta$ 
equals the corresponding number 
for $\Delta$ 
while 
the number of appearances of 1-vertices when traversing 
the interior face of $\tilde \Delta$ 
exceeds the corresponding number 
for $\Delta$ by 1. Therefore,
the equalities
$e(\tilde\Delta)=e(\Delta),$
$i(\tilde\Delta)=i(\Delta)+1$ hold.
Since by construction there are exactly $q_3+q_4+ ... +q_r$ stars of $\Delta$ for which 1-vertex is of 
valency 1 and $q_1-q_2\leq q_3 + ... +q_r$ by condition, after repeating this operation $q_1-q_2$ times 
we obtain a sunflower $\Sigma$ for which $q_i(\Sigma)=q_i,$ 
$1\leq i \leq r,$ and $i(\Sigma)=q_1.$ Notice that by construction $\Sigma$ has 
$q_2+ q_3 + ... +q_r-q_1$ 1-vertices of valency 1.

Now we are ready to construct the sunflower $\Omega$. First, observe that since
$e(\Sigma)=e(\Delta)=q_2+q_3 + ... +q_r$ in order to construct $\Omega$ for $s=q_2+q_3 + ... +q_r$
it is enough ``to turn inside out'' $\Sigma$ (see Fig. \ref{f5} where this operation is applied to the sunflower shown on Fig. \ref{f4}).
For $s,$ $1\leq s \leq q_2+q_3 + ... +q_r-1$
modify the sunflower $\Sigma$ as follows. Suppose first that $q_1< q_2+ q_3 + ... + q_r$. 
Then there exists a 1-vertex $u$ of $\Sigma$ of valency 2 such that 
the next $1$-vertex $v$, when traversing the exterior face of $\Sigma$ in the counter-clockwise direction, is of valency 1 (see Fig. \ref{f4}, where a  possible choice of $u$ and $v$ is shown). Indeed, consider an arbitrary 1-vertex $t$ of valency 2. If the condition above is not satisfied for $t$ then the next 1-vertex $t_1$ is also 
of valency 2. Check now the condition for $t_1$ and so on. 
Since the condition $q_1< q_2+ q_3 + ... + q_r$ 
implies that $\Sigma$ contains at least one 1-vertex of valency 1 continuing in this way we will arrive to the vertex needed  (recall that 1-vertices of valency 1 can not be adjacent to the interior face of $\Sigma$).

Now traverse the exterior face of $\Sigma$ in the counter-clockwise direction
starting from the vertex $v$  
till the moment when a $1$-vertex will appear for the $s$ time and denote this 1-vertex by $w.$ 
If the valency of $w$
is $2$ (see Fig. \ref{f4}, where $s=1$ and the corresponding vertex is denoted 
by $w_1$)  
then divide $w$ into two (not connected) $1$-vertices and glue one of them with $v$ as it shown on Fig. \ref{f6}  (note that if $s=q_2+q_3 + ... +q_r-1$ then $w=u$).
\begin{figure}[h]
\renewcommand{\captionlabeldelim}{.}
\medskip
\epsfxsize=12truecm
\centerline{\epsffile{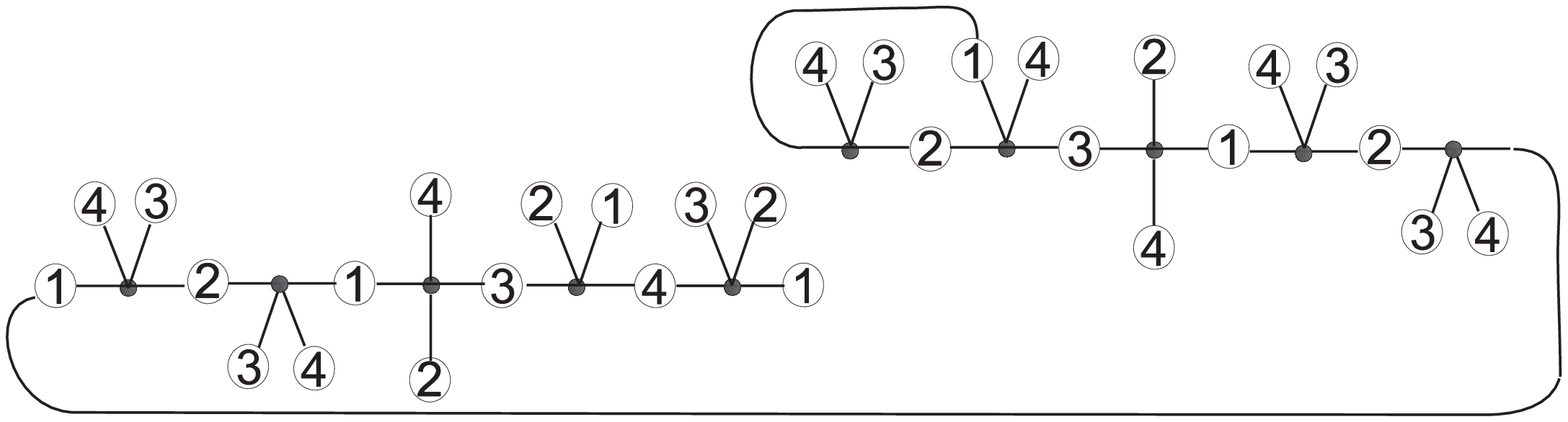}}
\caption{}\l{f6}
\medskip
\end{figure}

On the other hand, if the valency of $w$ is 1 (note that in this case necessarily 
$s<q_2+q_3 + ... +q_r-1$, see Fig. \ref{f4}, where $s=2$ and the corresponding vertex is denoted by $w_2$) then glue vertices $v$ and $w$ and then 
divide $u$ into two (not connected) $1$-vertices
as it is shown on Fig. \ref{f7}.
\begin{figure}[h]
\renewcommand{\captionlabeldelim}{.}
\medskip
\epsfxsize=12truecm
\centerline{\epsffile{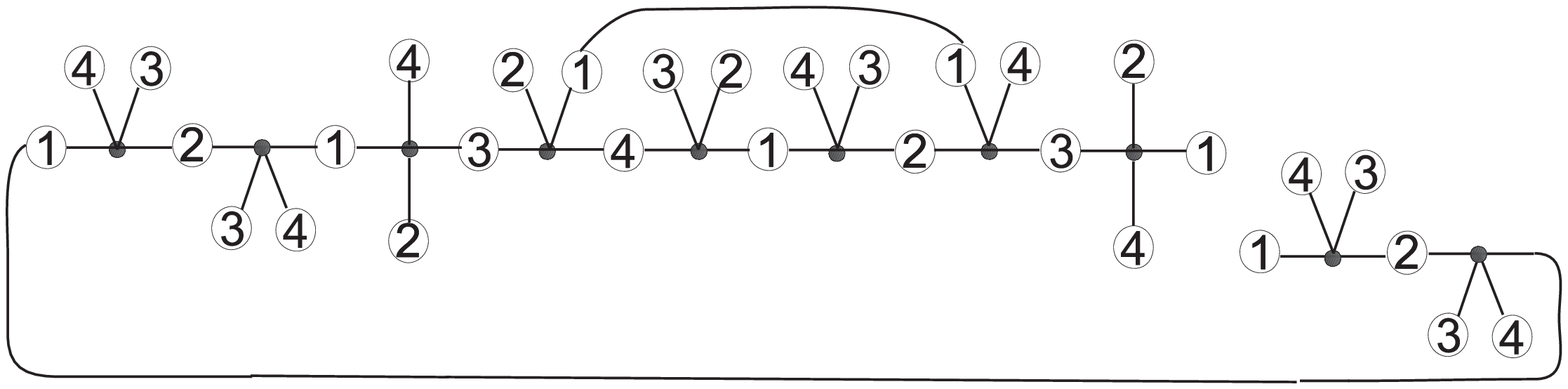}}
\caption{}\l{f7}
\medskip
\end{figure}
Clearly, in both cases we obtain a sunflower $\Omega$ for which $q_i(\Omega)=q_i,$ 
$1\leq i \leq r,$ and $i(\Omega)=s.$

To finish the proof we only must consider the case when $q_1= q_2+ q_3 + ... + q_r$
and $s$ satisfies
\be \l{blia} 1\leq s \leq q_2+ q_3 + ... + q_r-1.\ee 
Set $$\tilde q_1=q_2+ q_3 + ... + q_r-1, \ \   \tilde q_i=q_i, \ \ 
2\leq i \leq r.$$ 
Since $\tilde q_1<\tilde q_2+\tilde q_3 + ... +\tilde q_r$, 
for any number $s$ satisfying $1\leq s \leq \tilde q_2+\tilde q_3 + ... +\tilde q_r$
using the already proved
part of the proposition we can construct a sunflower  
$\tilde\Omega$ for which $q_i(\tilde\Omega)=\tilde q_i,$ $1\leq i \leq r,$ and $i(\tilde\Omega)=s.$
Furthermore, if $s$ satisfies 
\be \l{blia1} 1\leq s \leq \tilde q_2+\tilde q_3 + ... +\tilde q_r-1  \ee 
(that is if $\tilde\Omega$ is distinct from the sunflower shown on Fig. \ref{f5})
then by construction 
$\tilde\Omega$ contains a 1-vertex $y$ of valency 1 adjacent to the exterior face of $\tilde\Omega$ (see Fig. \ref{f6}, Fig. \ref{f7}).
Gluing now to the vertex $y$ a star we obtain a sunflower $\Omega$ 
for which $q_i(\Omega)=q_i,$ $1\leq i \leq r,$ and $i(\Omega)=s.$ Since the inequalities \eqref{blia}
and \eqref{blia1} are equivalent this proves the proposition. 
\begin{figure}[h]
\renewcommand{\captionlabeldelim}{.}
\medskip
\epsfxsize=12truecm
\centerline{\epsffile{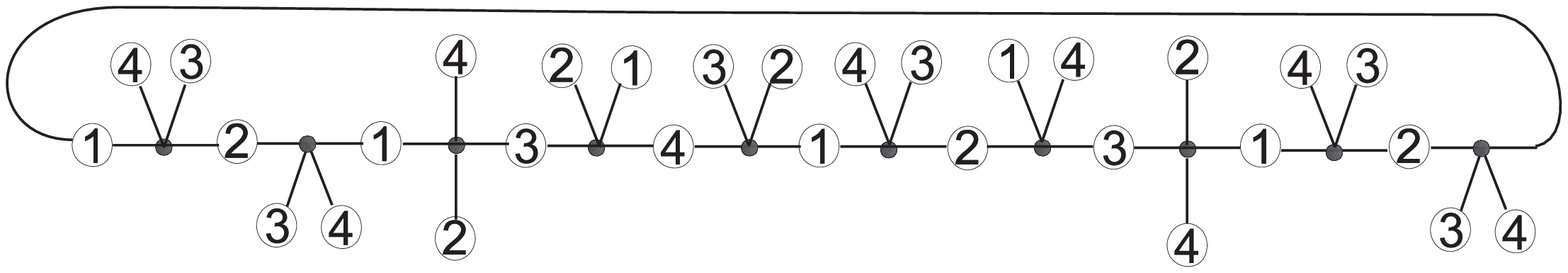}}
\caption{}\l{f5}
\medskip
\end{figure}

\bl \l{chain3} Any passport $\Pi$ for which
$s(\Pi) \leq  q_2(\Pi) + q_3(\Pi) + ... + q_r(\Pi)$ 
is realizable whenever $r(\Pi)>2.$ 

\el

\pr Suppose first that $q_1(\Pi) \leq  q_2(\Pi)+q_3(\Pi) + ... +q_r(\Pi).$ Then by proposition \ref{chain2}
there exists a sunflower $\Omega$ such that $i(\Omega)=s,$
$q_i(\Omega)=q_i(\Pi),$ $1\leq i \leq r,$ and all numerated vertices of $\Omega$ have valencies $\leq 2$.
Clearly, we can glue a number of stars to the vertices of valency 2 of $\Omega$
so that for the obtained constellation    
$\Omega_1$ to get \be \l{ssd} b_{i,j}(\Omega_1)=b_{i,j}(\Pi), \ \ \ \ 1\leq i \leq r, \ \ \ \ 1\leq j \leq q_i.\ee
Furthermore, since $\Omega$ is a sunflower we can glue the stars needed so that 
the constellation $\Omega_1$ also will be a sunflower
(see Fig. \ref{f8}, where $s(\Pi) <  q_2(\Pi) + q_3(\Pi) + ... + q_r(\Pi)$ and Fig. \ref{f11}, where 
$s(\Pi) = q_2(\Pi) + q_3(\Pi) + ... + q_r(\Pi)$). Then $i(\Omega_1)=s$ and
therefore the valency datum of $\Omega_1$ coincides with $\Pi$ 
(see the remark after lemma \ref{ed}). 
\begin{figure}[t]
\renewcommand{\captionlabeldelim}{.}
\medskip
\epsfxsize=12truecm
\centerline{\epsffile{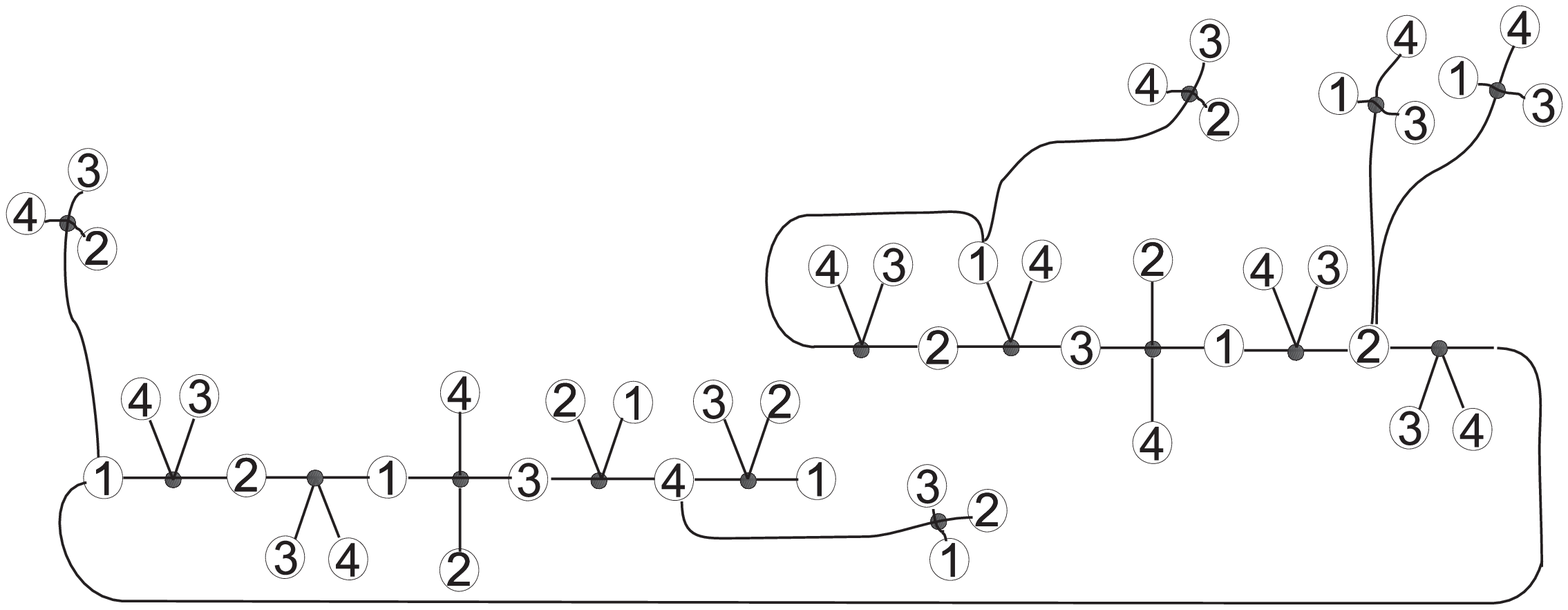}}
\caption{}\l{f8}
\medskip
\end{figure}

In the case when $q_1(\Pi) >  q_2(\Pi)+q_3(\Pi) + ... +q_r(\Pi)$ we act as follows.
In the beginning using proposition \ref{chain2} construct a sunflower $\Omega$ such that  $i(\Omega)=s$ and $$q_1(\Omega)=q_2(\Pi) + q_3(\Pi) + ... + q_r(\Pi), \ \ q_i(\Omega)=q_i(\Pi), \ \ 2\leq i \leq r.$$
Note that 
since $q_1(\Omega)=q_2(\Omega) + q_3(\Omega) + ... + q_r(\Omega)$ the construction of proposition \ref{chain2} implies that $\Omega$ contains no 
1-vertices of valency 1. 
In the next stage glue a number of stars to the vertices of valency 2 of $\Omega$ so that to obtain a sunflower $\Omega_1$ for which $i(\Omega_1)=s$ and 
$$b_{i,j}(\Omega_1)=b_{i,j}(\Pi), \ \ \ \ 2\leq i \leq r, \ \ \ \ 1\leq j \leq q_i,$$
while
$$\{b_{1,1}(\Omega_1), b_{1,2}(\Omega_1),\, ... \, , b_{1,q_1(\Omega_1)}(\Omega_1)\}=\{b_{1,l+1}(\Pi), b_{1,l+2}(\Pi),\, ...\, , b_{1,q_1(\Pi)}(\Pi)\},$$
where $$l=q_1(\Pi)-(q_2(\Pi)+q_3(\Pi)+ ... + q_r(\Pi))$$ (see Fig. \ref{f9}, where $s(\Pi) <  q_2(\Pi) + q_3(\Pi) + ... + q_r(\Pi)$).
\begin{figure}[h]
\renewcommand{\captionlabeldelim}{.}
\medskip
\epsfxsize=8truecm
\centerline{\epsffile{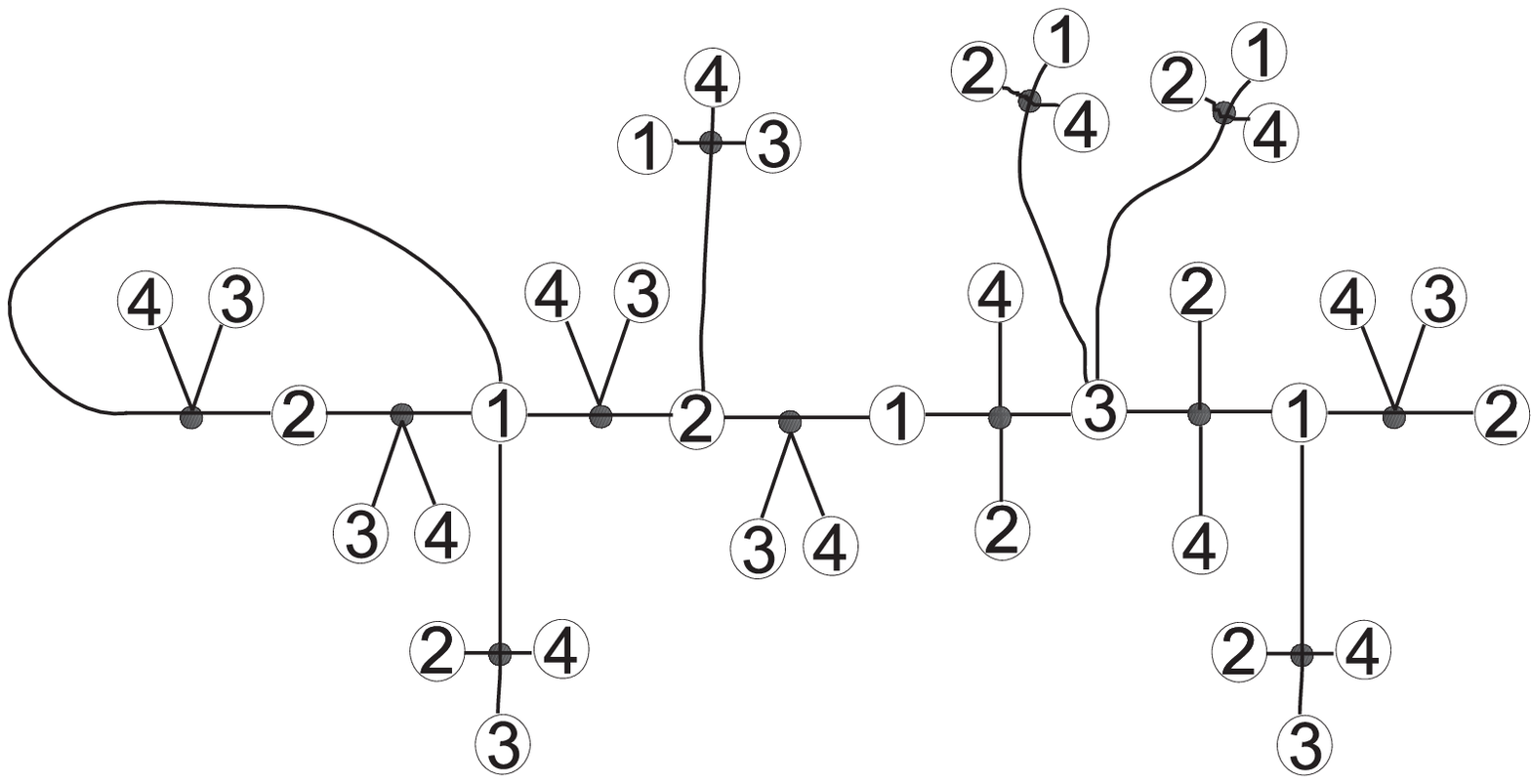}}
\caption{}\l{f9}
\medskip
\end{figure}

Since $\Omega$ have no 1-vertices of valency 1 it is easy to see that for the number $\nu$ of 1-vertices of valency 1 of $\Omega_1$ 
the equality
$$\nu=\sum_{i=2}^r \sum_{j=1}^{q_i(\Pi)} (b_{i,j}(\Pi)-2)$$ holds. Note that all these 1-vertices are adjacent 
to the exterior face of $\Omega_1.$ 
Since lemma \ref{gop} implies that $\nu \geq l$ as the last stage of our construction we can  
glue $\nu-l$ stars to the 1-vertices of valency 1 of $\Omega_1$ so that to obtain 
a sunflower $\Omega_2$ for which $i(\Omega_2)=s,$ 
$q_i(\Omega_2)=q_i(\Pi),$ $1\leq i \leq r,$ and 
$$b_{i,j}(\Omega_2)=b_{i,j}(\Pi), \ \ \ \ 1\leq i \leq r, \ \ \ \ 1\leq j \leq q_i$$ 
(see Fig. \ref{f10}, where $s(\Pi) <  q_2(\Pi) + q_3(\Pi) + ... + q_r(\Pi)$ and Fig. \ref{f13}, where  
$s(\Pi) = q_2(\Pi) + q_3(\Pi) + ... + q_r(\Pi)$).
\begin{figure}[h]
\renewcommand{\captionlabeldelim}{.}
\medskip
\epsfxsize=8truecm
\centerline{\epsffile{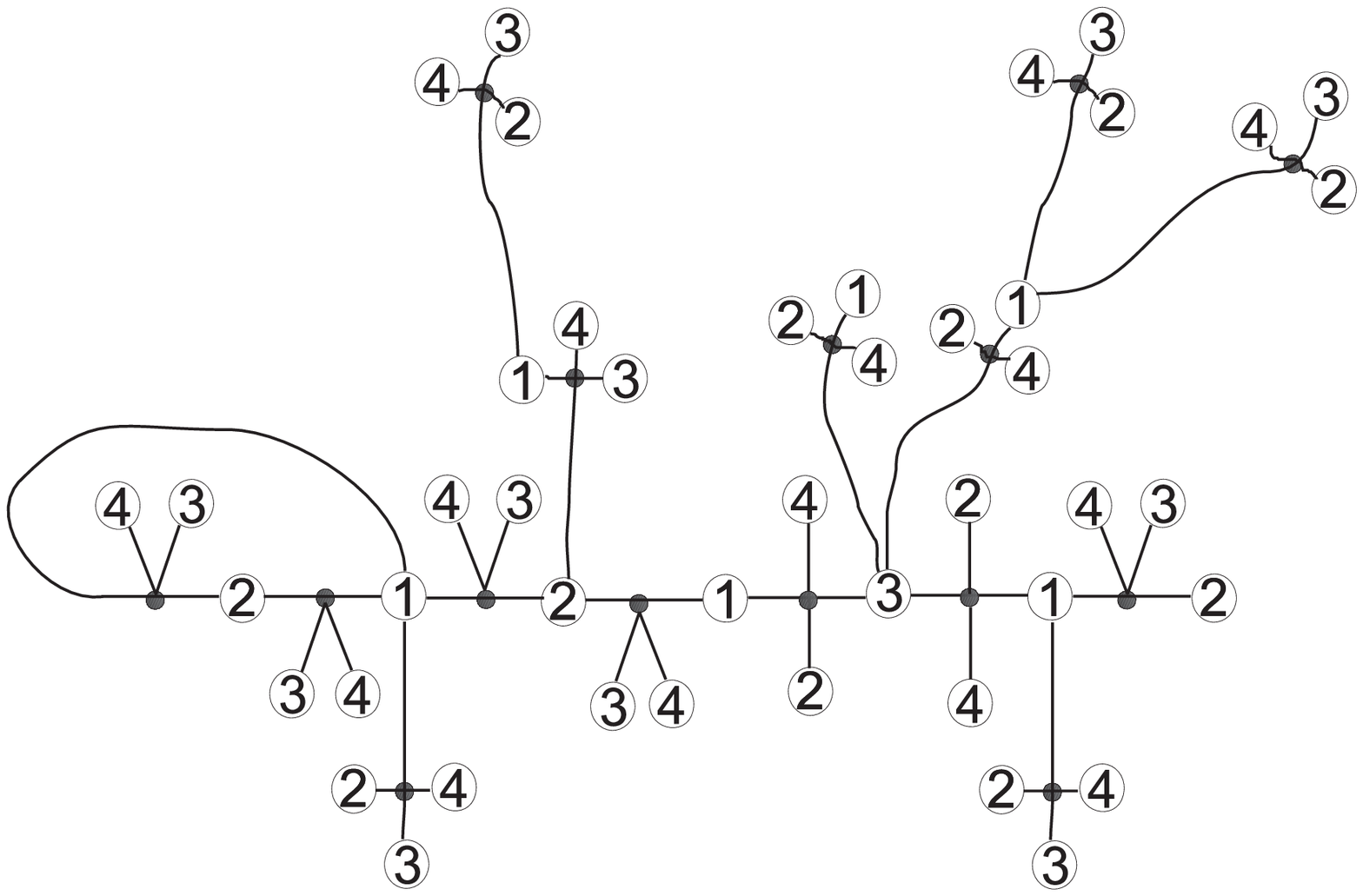}}
\caption{}\l{f10}
\medskip
\end{figure}

\bp \l{4} Any passport $\Pi$ for which 
$q_1(\Pi)\leq q_2(\Pi) + q_3(\Pi) + ... + q_r(\Pi)$ 
is realizable whenever $r(\Pi)> 2.$ 
\ep

\pr 
In view of lemma \ref{chain3} we only must consider the case when $s$ satisfies
$$q_2(\Pi)+q_3(\Pi) + ... +q_r(\Pi) < s \leq n/2.$$
Let $\Omega$ be a sunflower 
such that $\Omega_i=\Pi_i$, $1\leq i \leq r,$ and 
$$i(\Omega)=q_2(\Pi)+q_3(\Pi) + ... +q_r(\Pi)$$ constructed in lemma \ref{chain3}.
Since $q_1(\Pi)\leq q_2(\Pi) + q_3(\Pi) + ... + q_r(\Pi)$, $\Omega$ has the form
shown on Fig. \ref{f11}.

Observe that is we ``shift'' any of branches of $\Omega$ from outside to inside (see Fig. \ref{f12}) 
then we obtain a constellation $\tilde\Omega$ with $q_i(\tilde\Omega)=q_i(\Pi),$ $1\leq i \leq r,$ and
$$i(\tilde\Omega)=q_2(\Pi)+q_3(\Pi) + ... +q_r(\Pi)+1.$$
It is clear that repeating this operation
we can obtain a constellation $\Omega_1$ with $q_i(\Omega_1)=q_i(\Pi),$ $1\leq i \leq r,$ and $i(\Omega_1)$ equal to any $s$ which satisfies 
$$q_2(\Pi)+q_3(\Pi) + ... +q_r(\Pi)+1 \leq s \leq \mu,$$ where
\be \l{fo1} \mu=q_2(\Pi)+q_3(\Pi) + ... +q_r(\Pi)+
\sum_{i=1}^r \sum_{j=1}^{q_i(\Pi)} (b_{i,j}(\Pi)-2).\ee
\begin{figure}[t]
\renewcommand{\captionlabeldelim}{.}
\medskip
\epsfxsize=12truecm
\centerline{\epsffile{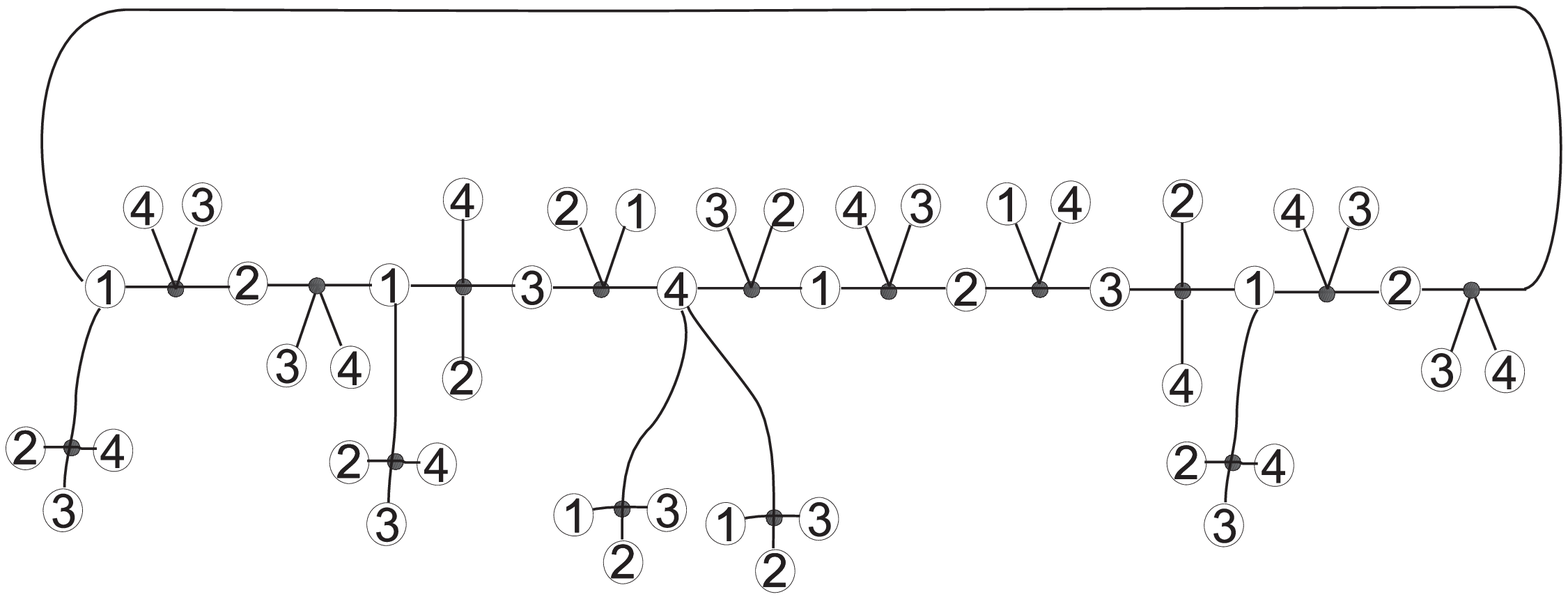}}
\caption{}\l{f11}
\medskip
\end{figure}

So, 
to finish the proof we only must show that $\mu \geq n/2.$ 
Since by lemma \ref{gop} $$\mu=
\sum_{j=1}^{q_1(\Pi)} (b_{1,j}(\Pi)-2)+ e_1(\Pi)+q_1(\Pi)=n-e_1(\Pi)-2q_1(\Pi)+e_1(\Pi)+q_1(\Pi)=$$ \be \l{ffo2} =n-q_1(\Pi),\ee
it follows from the obvious equality $q_1(\Pi)\leq n/2$  
that \be \l{fo2} \mu \geq n/2.\ee

\begin{figure}[h]
\renewcommand{\captionlabeldelim}{.}
\medskip
\epsfxsize=12truecm
\centerline{\epsffile{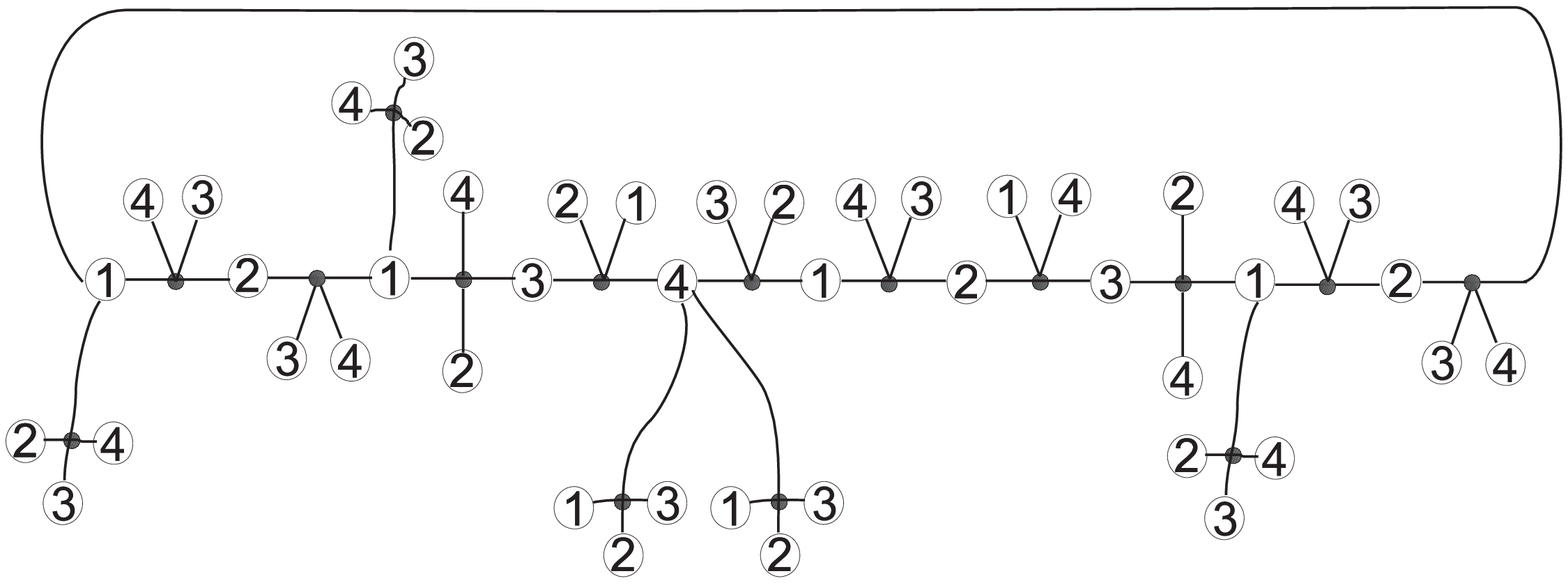}}
\caption{}\l{f12}
\medskip
\end{figure}


\bl \l{cal} Let $u_1, u_2, ... u_l,t$ be integers such that 
$1<u_1\leq u_2\leq ... \leq u_l$ and $t\geq 1.$ Then the equation 
\be \l{ur}
s=y+x_1u_1+x_2u_2+ ... + x_lu_l, \ee has a solution in $x_i,y$ 
with $x_i \in \{0,1\},$ $1\leq i \leq l,$ and $y \in \{0,1,2, ...\, ,t\}$ for any 
$s$ satisfying $0\leq s \leq t+u_1+u_2+ ... +u_l$ if and only if  
for any $k,$ $1\leq k\leq l,$ the inequality \be \l{us} t+\sum_{i=1}^{k-1} u_i \geq u_k-1 \ee holds. In particular, the condition $ t \geq u_l-1$ is sufficient.

\el

\pr First notice that condition \eqref{us} is necessary since if \eqref{us} fails 
to be true say for $k=h$ 
then equation \eqref{ur} has no solutions for $s=u_h-1.$ Indeed, 
since $s< u_h\leq u_{h+1}\leq ... \leq u_l$ if such a solution exists then necessary 
$x_i=0$ for $i\geq h.$ On the other hand, since $t+\sum_{i=1}^{h-1} u_i < u_h-1,$
the inequality $y+x_1u_1+x_2u_2+ ... + x_{h-1}u_{h-1}< s$ holds
for any choice of $x_i,$ $1\leq i \leq h-1,$ and $y,$ $0\leq y \leq t.$

To prove the sufficiency of \eqref{us} we use the induction by $l.$ For $l=1$ the lemma is obvious. Suppose that it holds for $l=n$ and prove it for $l=n+1.$  
If $s$ satisfies $0\leq s \leq t+u_1+u_2+ ... +u_n$ then the statement is true since by the inductive hypothesis
there exist $x_i,$ $1\leq i \leq n,$ and $y,$ $0\leq y \leq t,$ such that
$$s=y+x_1u_1+x_2u_2+ ... + x_{n}u_{n}.$$ 

On the other hand, if
\be \l{ner} t+u_1+u_2+ ... +u_n < s \leq t+u_1+u_2+ ... +u_n +u_{n+1}\ee
then \eqref{us} taken for $k=l=n+1$ implies that $s=u_{n+1}+\tilde s$ for some $\tilde s \geq 0.$
Furthermore, since \eqref{ner} implies that $0 \leq \tilde s \leq t+u_1+u_2+ ... +u_n,$  
the inductive hypothesis implies that there exist $x_i,$ $1\leq i \leq n,$ and $y,$ $0\leq y \leq t,$ such that
$$\tilde s=y+x_1u_1+x_2u_2+ ... + x_{n}u_{n}$$ and hence 
$$s=y+x_1u_1+x_2u_2+ ... + x_{n}u_{n}+x_{n+1}u_{n+1}$$
with $x_{n+1}=1.$

\bl \l{ss} Any  passport $\Pi$ for which
$q_1(\Pi) > q_2(\Pi) + q_3(\Pi) + ... + q_r(\Pi)$
and $s$ satisfies $q_2(\Pi) + q_3(\Pi) + ... + q_r(\Pi)< s(\Pi) \leq n/2$
is realizable whenever $\Pi_1\neq \{2,2, ... ,2\}$ and $r(\Pi)>2.$
\el

\pr 
Let $\Omega$ be a sunflower such that $\Omega_i=\Pi_i,$ $1\leq i \leq r,$ and $i(\Omega)=q_2(\Pi)+q_3(\Pi) + ... +q_r(\Pi)$ constructed in lemma \ref{chain3}. 
In view of the inequality  
$q_1(\Pi) > q_2(\Pi) + q_3(\Pi) + ... + q_r(\Pi)$, $\Omega$ 
has the form shown on Fig. \ref{f13} 
\begin{figure}[t]
\renewcommand{\captionlabeldelim}{.}
\medskip
\epsfxsize=11truecm
\centerline{\epsffile{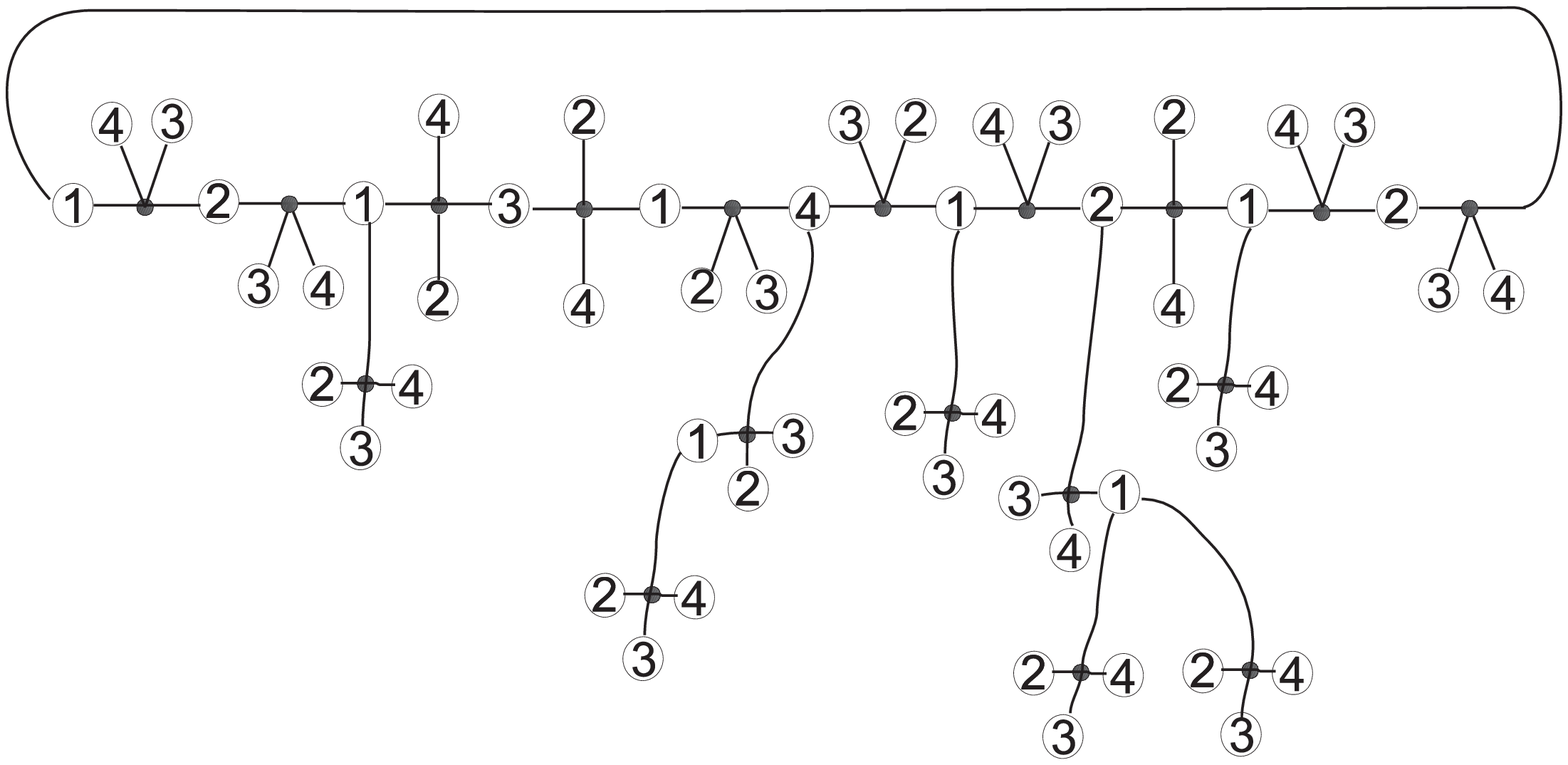}}
\caption{}\l{f13}
\medskip
\end{figure}
and by construction admits two types of branches. First $\Omega$ has $$l=q_1(\Pi)-(q_2(\Pi) + q_3(\Pi) + ... + q_r(\Pi))$$ ``long''
branches $\lambda_i$ 
for which
$\vert \lambda_i \vert= b_{1,i},$ $1\leq i \leq l.$ 
Second $\Omega$ has \be \l{ff} t= \sum_{j=l+1}^{q_1(\Pi)}(b_{1,j}(\Pi)-2) +\sum_{i=2}^r\sum_{j=1}^{q_i(\Pi)}(b_{i,j}(\Pi)-2)-l \ee
``short'' branches $\mu_j,$ $1\leq j \leq t,$ for which $\vert \mu_j \vert=1.$ 
Note that in view of lemma \ref{gop} we have: $$ t= \sum_{j=l+1}^{q_1(\Pi)}(b_{1,j}(\Pi)-2)+e_1(\Pi).$$

Clearly, shifting a number of branches $\lambda_i,$ $1\leq i \leq l,$ $\mu_j,$ $1\leq j \leq t,$
from outside to inside (see Fig. \ref{f14})  
we can obtain 
a constellation $\Omega_1$ such that ${\Omega_1}_i=\Pi_i$, $1\leq i \leq r,$ and $i(\Omega_1)=s,$
where $s$ is any number which can be represented as the sum
\be \l{sss} s=q_2(\Pi)+q_3(\Pi) + ... +q_r(\Pi)+y+x_1b_{1,1}(\Pi)+x_2b_{1,2}(\Pi)+ ... + x_lb_{1,l}(\Pi) \ee for some $x_i,y$ 
with $x_i \in \{0,1\},$ $1\leq i \leq l,$ and $y \in \{0,1,2, ...\, ,t\}.$
\begin{figure}[t]
\renewcommand{\captionlabeldelim}{.}
\medskip
\epsfxsize=11truecm
\centerline{\epsffile{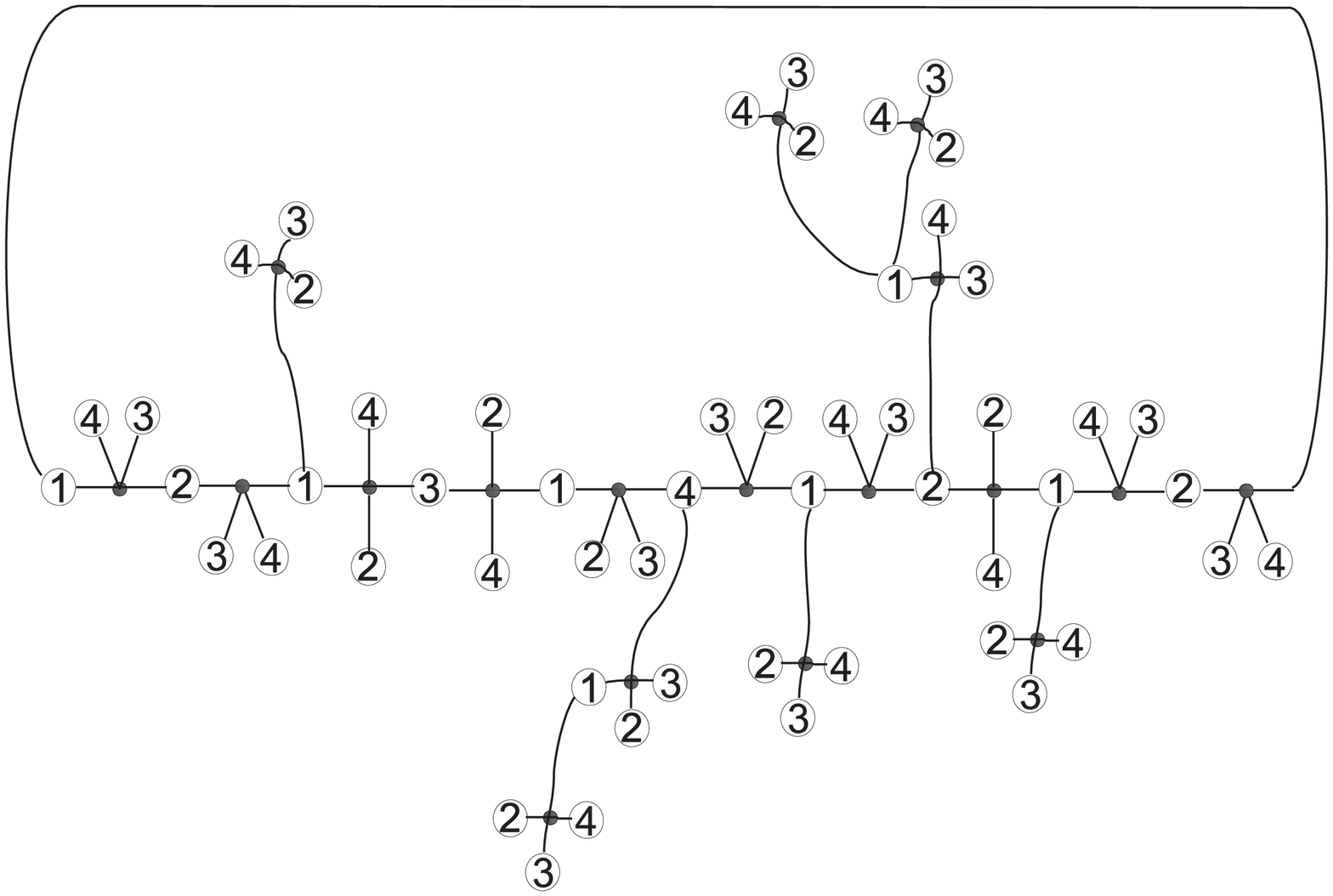}}
\caption{}\l{f14}
\medskip
\end{figure}
Furthermore, since for the maximal possible value $s_{\rm max}$ of $s$ we have:
$$s_{\rm max}=q_2(\Pi)+q_3(\Pi) + ... +q_r(\Pi)+\sum_{j=1}^{l}b_{1,j}(\Pi) +t=$$ $$=q_1(\Pi)-l+\sum_{j=1}^{l}b_{1,j}(\Pi)+\sum_{j=l+1}^{q_1(\Pi)}(b_{1,j}(\Pi)-2)+e_1(\Pi)=$$$$
=q_1(\Pi)-l+\sum_{j=1}^{q_1(\Pi)}b_{1,j}(\Pi)-2(q_1(\Pi)-l)+e_1(\Pi)=$$$$ =\sum_{j=1}^{q_1(\Pi)}b_{1,j}(\Pi)+e_1(\Pi)
-(q_1(\Pi)-l)=n-(q_2(\Pi)+q_3(\Pi) + ... +q_r(\Pi)),$$ it follows from $$q_2(\Pi)+q_3(\Pi) + ... +q_r(\Pi)<q_1(\Pi)\leq n/2$$ that \be \l{fo4} s_{\rm max}>n/2. \ee

Therefore, in order to prove the lemma we only must show that $s$ can take any value
between $0$ and $s_{\rm max}.$
By lemma \ref{cal} it is enough to establish
that \be \l{los} t=\sum_{j=l+1}^{q_1(\Pi)}(b_{1,j}(\Pi)-2)+e_1(\Pi)\geq b_{1,l}(\Pi)-1.\ee
Since the condition $r>2$ implies that \be \l{fo5} q_1(\Pi)-l=q_2(\Pi) + q_3(\Pi) + ... + q_r(\Pi)\geq 2 \ee we have: 
$$\sum_{j=l+1}^{q_1(\Pi)}(b_{1,j}(\Pi)-2)+e_1(\Pi)\geq b_{1,q_1}(\Pi)+ b_{1,q_1-1}(\Pi)-4+e_1(\Pi).$$
Furthermore, since $\Pi_1\neq \{2,2, ... ,2\}$ at least one of the inequalities $b_{1,q_1}(\Pi)\geq 3,$ $e_1(\Pi)>1$ holds.
In both cases we have: \be \l{loss} b_{1,q_1}(\Pi)+ b_{1,q_1-1}(\Pi)-4+e_1(\Pi)\geq b_{1,q_1-1}(\Pi)-1 \geq b_{1,l}(\Pi)-1.\ee

\bp \l{5} Any passport $\Pi$ for which $q_1(\Pi)> q_2(\Pi) + q_3(\Pi) + ... + q_r(\Pi)$
is realizable whenever $r> 2.$ 
\ep 

\pr If $s\leq q_2(\Pi)+q_3(\Pi) + ... +q_r(\Pi)$ then the proposition follows from lemma 
\ref{chain3}. If $q_2(\Pi) + q_3(\Pi) + ... + q_r(\Pi)< s \leq n/2$ and $\Pi_1\neq \{2,2, ... ,2\}$ then the proposition follows from lemma \ref{ss}. Therefore we only must consider the case when $$\Pi_1= \{2,2, ... ,2\}, \ \ \ \ q_2(\Pi) + q_3(\Pi) + ... + q_r(\Pi)< s \leq n/2.$$

Let $\Omega$ be a sunflower 
such that $\Pi_i(\Omega)=\Pi_i,$ $1\leq i \leq r,$ and $i(\Omega)=q_2(\Pi) + q_3(\Pi) + ... + q_r(\Pi)$ 
constructed in lemma \ref{chain3} (see Fig. \ref{f15}). Since $\Pi_1= \{2,2, ... ,2\}$ it has a 
more restrictive form than the one shown on Fig. \ref{f13} in particular the branches of $\Omega$ can grow only from non 1-vertices and have weight 2.   
\begin{figure}[h]
\renewcommand{\captionlabeldelim}{.}
\medskip
\epsfxsize=10truecm
\centerline{\epsffile{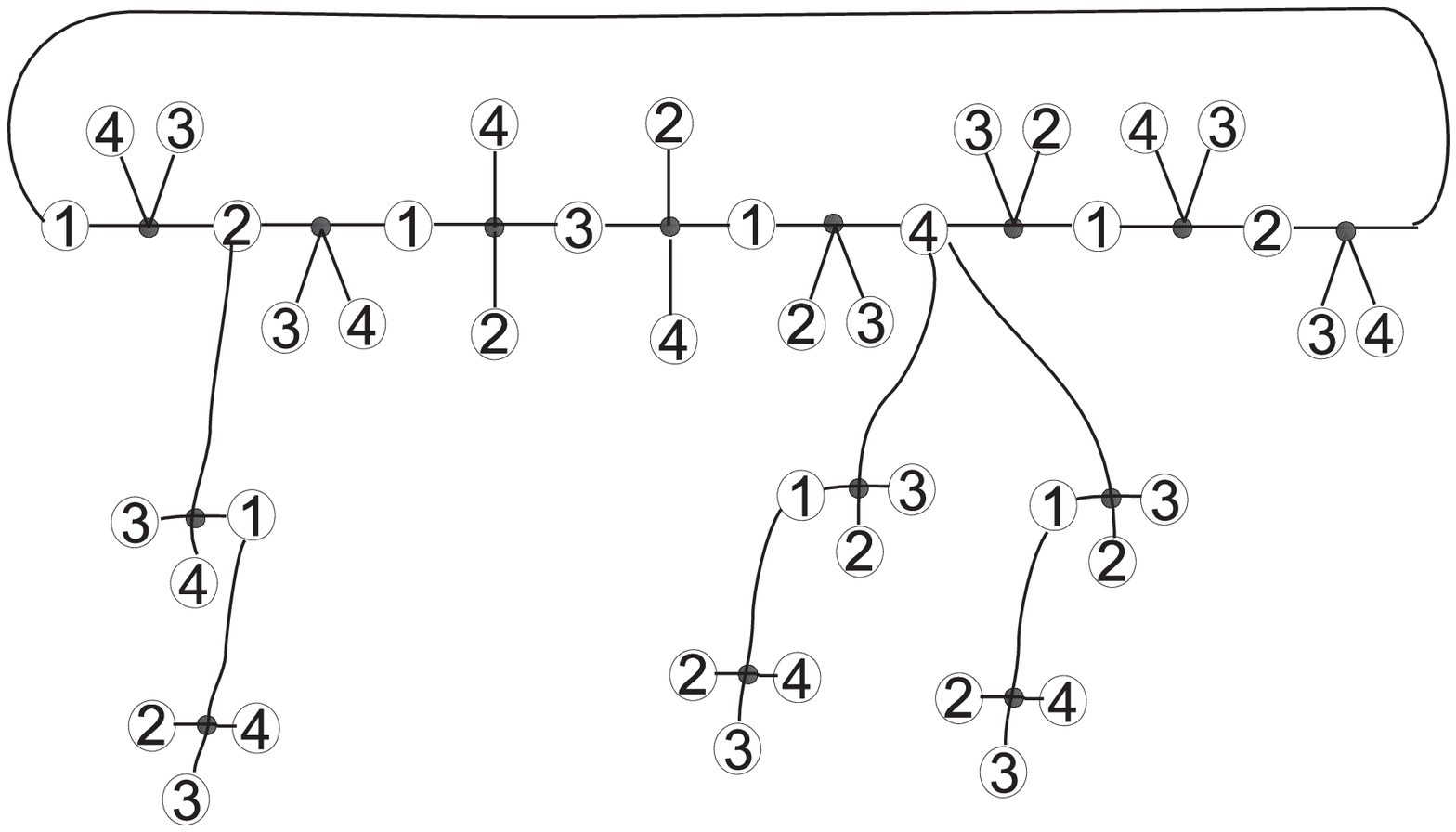}}
\caption{}\l{f15}
\medskip
\end{figure}
As above shifting these branches from outside to 
inside  we can obtain 
a constellation $\Omega_1$ such that $\Pi_i(\Omega_1)=\Pi_i,$ $1\leq i \leq r,$ and
$i(\Omega_1)=s,$ where $s$ is any number which has the form
$$s= q_2(\Pi)+q_3(\Pi) + ... +q_r(\Pi)+2k,\ \ \ \ 0\leq k \leq \sum_{i=2}^r\sum_{j=1}^{q_i(\Pi)}(b_{i,j}(\Pi)-2).$$
Since in view of lemma \ref{gop} for the maximal possible value $s_{\rm max}$ of $s$ we have:
$$s_{\rm max}=2q_1(\Pi) - (q_2(\Pi)+q_3(\Pi) + ... +q_r(\Pi))$$
and $q_1(\Pi) > q_2(\Pi)+q_3(\Pi) + ... +q_r(\Pi)$
the inequality $s_{\rm max}>q_1(\Pi)$ holds.
It follows now from $q_1(\Pi)=n/2$ that
$ s_{\rm max}> n/2$ and therefore $\Pi$ is realizable whenever 
$s= q_2(\Pi)+q_3(\Pi) + ... +q_r(\Pi)\ (\mod 2).$

In order to treat the case when $$s=1+ q_2(\Pi)+q_3(\Pi) + ... +q_r(\Pi)\ (\mod 2), \ \ \ \ q_2(\Pi) + q_3(\Pi) + ... + q_r(\Pi)< s \leq n/2$$ we
act as follows. In the beginning using the already proved part of the proposition 
construct a sunflower $\Omega_2$ such that $\Pi_i(\Omega_2)=\Pi_i,$ $1\leq i \leq r$ and
$i(\Omega_2)=s-1.$ 
Recall that by construction (see proposition \ref{chain2}) the cycle $\c(\Omega_2)$ of $\Omega_2$ possesses the following property: the 1- and non 1-vertices of $\c(\Omega_2)$ alternate and if a non 1-vertex $v$ follows a non 1-vertex $u$,
when traversing $\c(\Omega_2)$ in the counter-clockwise direction, then the number of $v$ is greater than
the number of $u$ unless $v$ is a 2-vertex
(see Fig. \ref{rys}). In particular, since $r>2$ we can find a
pair of vertices $u,v$ such that $u$ is a 2-vertex, $v$ is a 3-vertex, and $v$ follows $u.$ 

\begin{figure}[h]
\renewcommand{\captionlabeldelim}{.}
\medskip
\epsfxsize=10truecm
\centerline{\epsffile{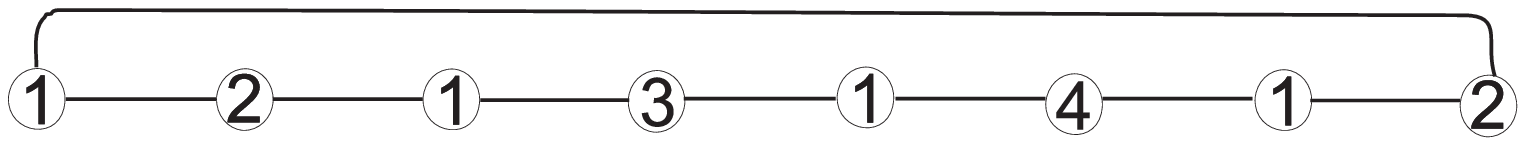}}
\caption{}\l{rys}
\medskip
\end{figure}

Consider the corresponding adjacent stars $S,R$ of $\Omega_2$ (see Fig. \ref{f16},a) and perform the following operation: remove $S,R$ and glue instead two new stars shown on Fig. \ref{f16},b
leaving the branches possibly growing from $u$ and $v$ (denoted by dotted lines) unchanged (see Fig. \ref{f17}, where this operation is applied to the constellation shown on 
Fig. \ref{f15}).
\begin{figure}[h]
\renewcommand{\captionlabeldelim}{.}
\medskip
\epsfxsize=11truecm
\centerline{\epsffile{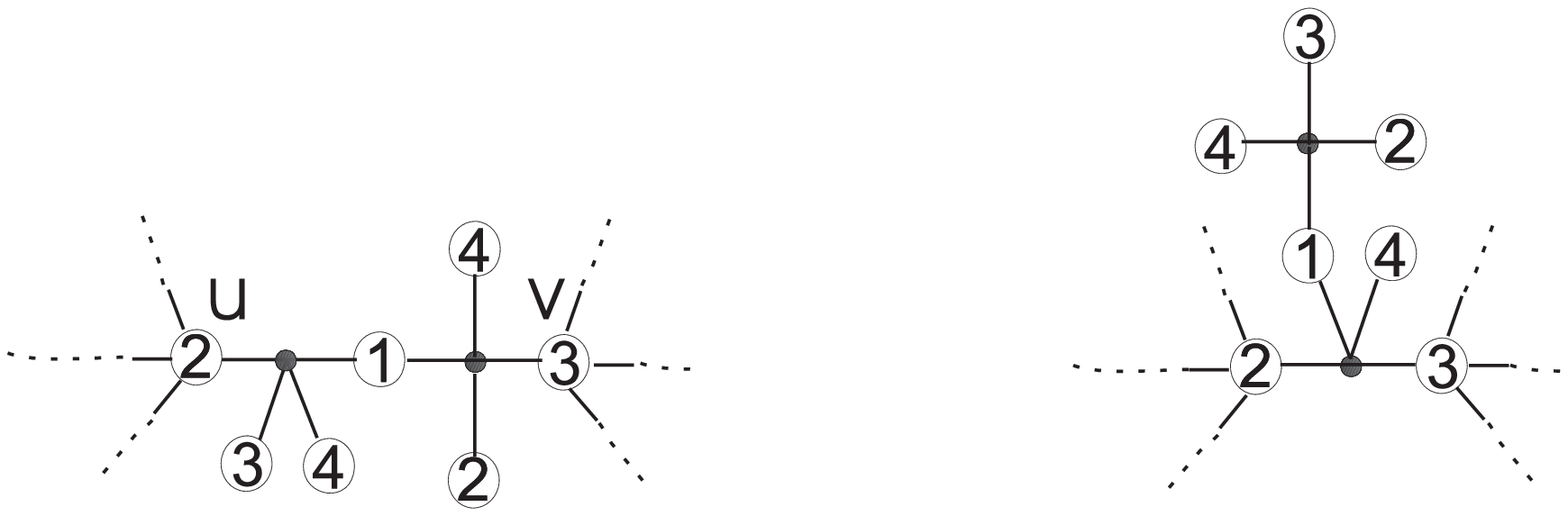}}
\caption{}\l{f16}
\medskip
\end{figure}
 
Taking into account that the branches of $\Omega_2$ 
can grow only from the non 1-vertices of valency $2$ it is easy to see that this operation is well defined and that as a result we obtain a constellation    
$\Omega_3$ for which ${\Omega_3}_i=\Pi_i,$ $1\leq i \leq r,$ and $i(\Omega_3)=s$.

\begin{figure}[t]
\renewcommand{\captionlabeldelim}{.}
\medskip
\epsfxsize=9truecm
\centerline{\epsffile{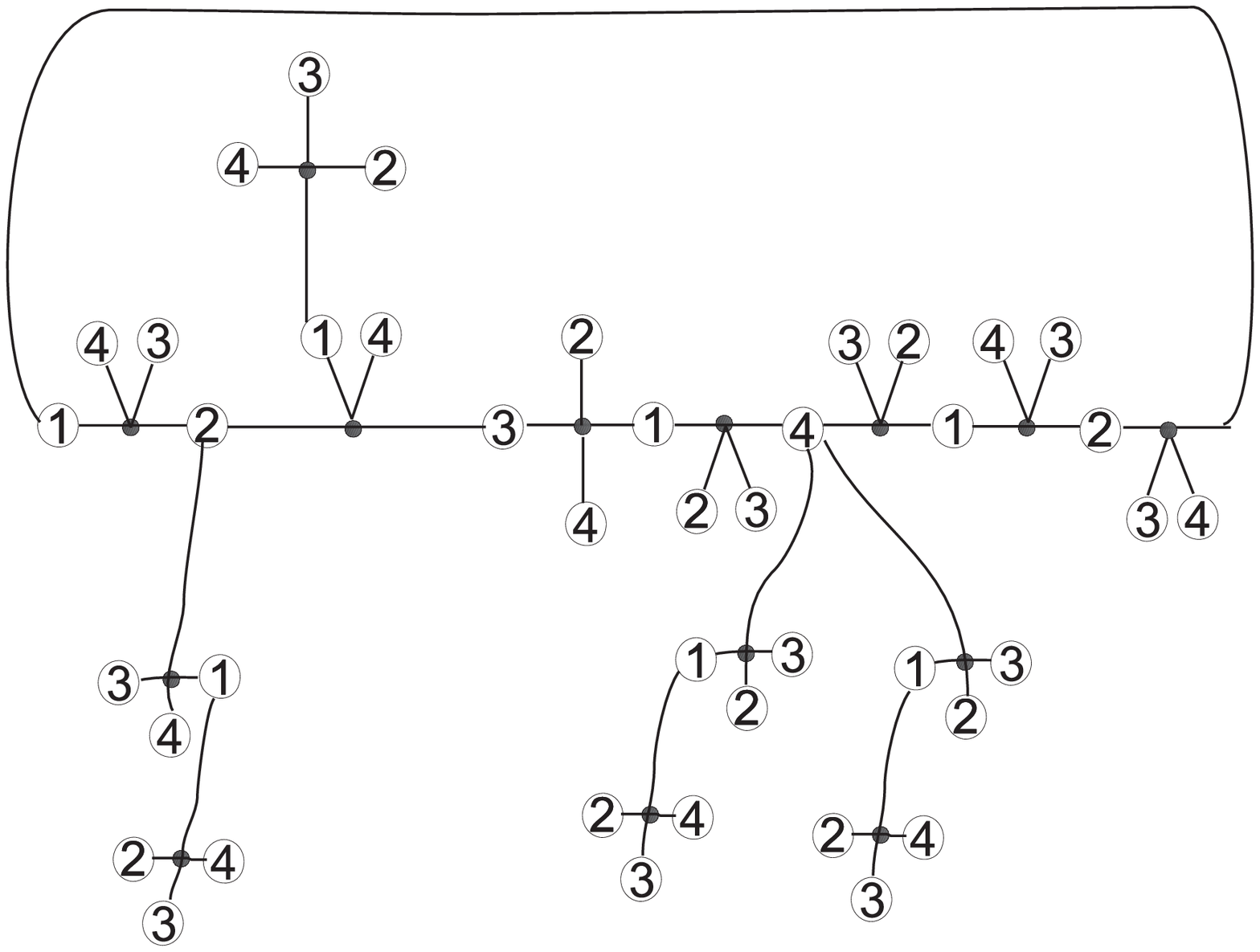}}
\caption{}\l{f17}
\medskip
\end{figure}

\bt Any passport for which $r(\Pi)>2$ is realizable.
\et

\pr Follows from propositions \ref{4}, \ref{5}.

\section{Passports with $r=2$.} 
In this section we will picture all constellations in the form 
of bicolored graphs 
(see subsection 2.1). 
\bl \l{ll} Let $\Pi$ be a passport such that $r(\Pi)=2$ and either 
$\Pi_1\neq\{2,2,\, ...\, ,2\}$
or $\Pi_2\neq \{2,2,\, ...\, ,2\}.$ Then either $b_{1,q_1}(\Pi)>2$ or $b_{2,q_2}(\Pi)>2.$
\el

\pr If $\Pi_1\neq\{2,2,\, ...\, ,2\}$ then either $b_{1,q_1}(\Pi)>2$ or $e_1(\Pi) >0$. 
On the other hand, by lemma \ref{gop}
\be\l{eche}\sum_{j=1}^{q_2(\Pi)}(b_{2,j}(\Pi)-2)=e_1(\Pi)+q_1(\Pi)-q_2(\Pi).\ee
Since it is assumed that $q_1(\Pi)\geq q_2(\Pi)$ it follows that
if $e_1(\Pi) >0$ then $b_{2,q_2}(\Pi)>2$. 

If $\Pi_1=\{2,2,\, ...\, ,2\}$
then \eqref{rh1} implies that $e_2(\Pi)+q_2(\Pi)=n/2.$ 
Furthermore, since $\Pi_2\neq\{2,2,\, ...\, ,2\}$ 
the inequality $q_2(\Pi)<n/2$ holds. Therefore,   
$e_2(\Pi)>0$ and hence $b_{2,q_2}(\Pi)>2$ since otherwise  $$\sum_{j=1}^{p_2(\Pi)}a_{2,j}=e_2(\Pi)+2q_2(\Pi)=e_2(\Pi)+2(n/2-e_2(\Pi))=n-e_2(\Pi)<n.$$


\bl \l{chai} Any passport $\Pi$ with $r(\Pi)=2$ for which 
$s \leq   q_2(\Pi)$ 
is realizable whenever $\Pi_1\neq\{2,2,\, ...\, ,2\}$ or $\Pi_2\neq \{2,2,\, ...\, ,2\}.$
\el

\pr 
To construct the constellation needed we act similarly to the case when $r>2$
with some simplifications. Suppose first that $s < q_2(\Pi).$ In the beginning construct a sunflower $\Omega$, which has one vertex of valency 1, one 
vertex of valency 3, and all other vertices of valency 2, 
such that $q_1(\Omega)=q_2(\Omega)=q_2(\Pi)$ and $i(\Omega)=s$ as it shown on Fig. \ref{f18}
(the number of vertex of valency 3 coincides with $i=1,2$ for which 
$b_{i,q_i}(\Pi)>2$).
\begin{figure}[h]
\renewcommand{\captionlabeldelim}{.}
\medskip
\epsfxsize=2truecm
\centerline{\epsffile{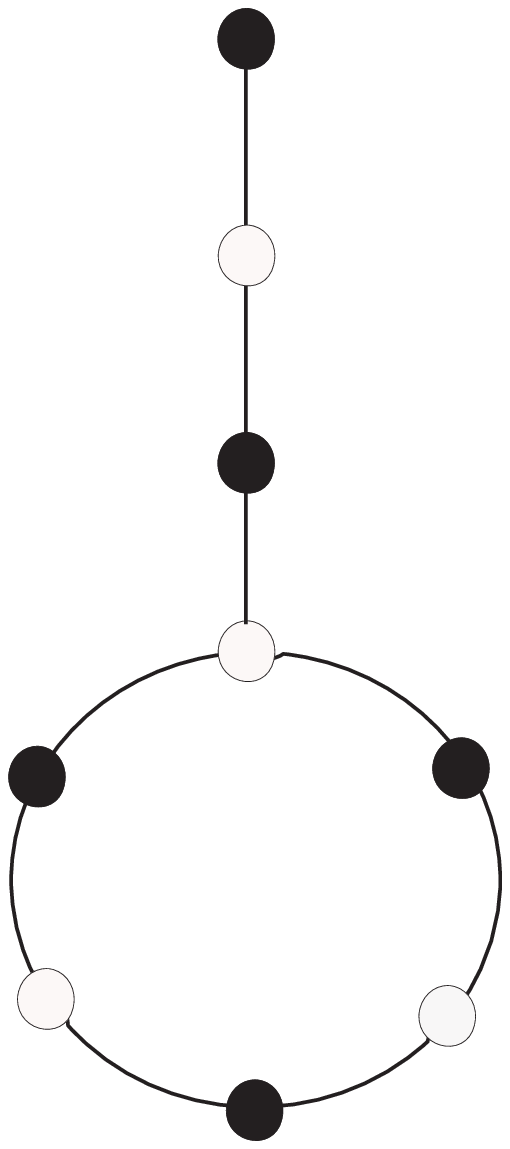}}
\caption{}\l{f18}
\medskip
\end{figure}

If $q_1(\Pi) = q_2(\Pi)$ then in order to construct 
a sunflower    
$\Sigma$ for which 
${\Sigma}_i=\Pi_i,$ $1\leq i \leq r,$ and
$i(\Sigma) = s$ 
it is enough 
to glue a number of edges to the vertices of valency 2 and 3 of the sunflower $\Omega$ 
(see Fig. \ref{f19}).
\begin{figure}[h]
\renewcommand{\captionlabeldelim}{.}
\medskip
\epsfxsize=3truecm
\centerline{\epsffile{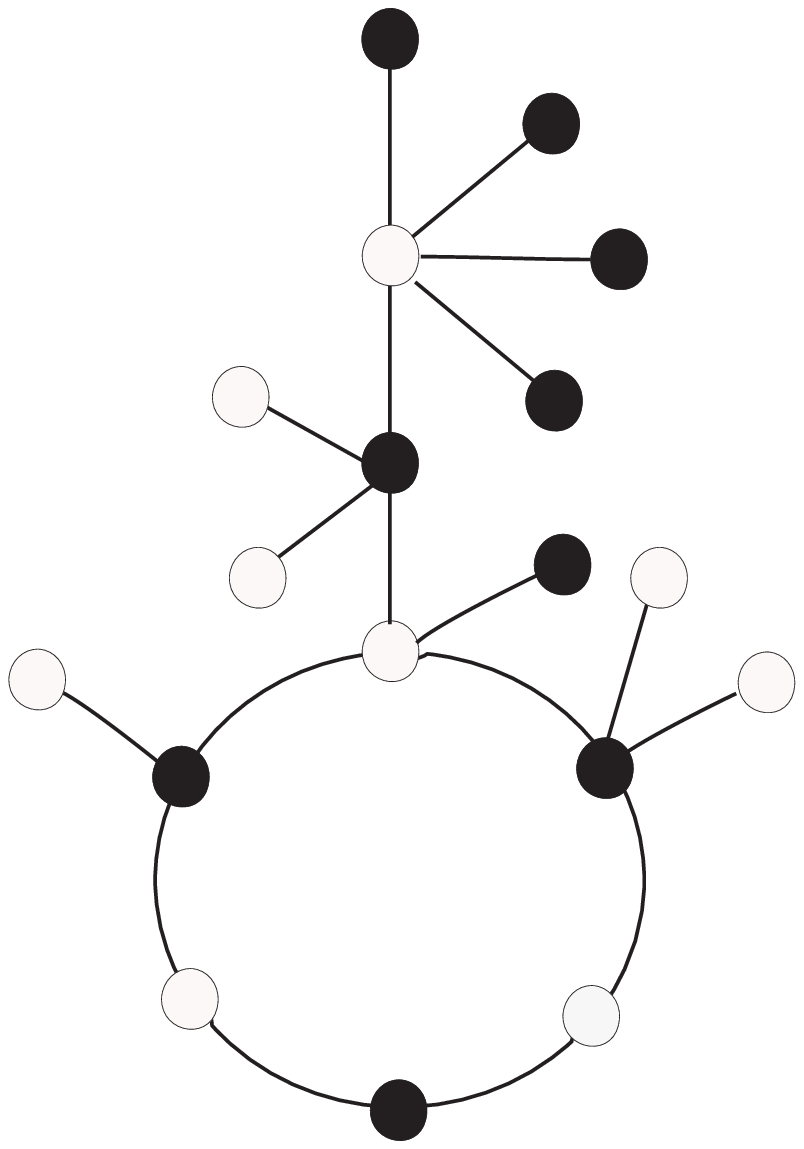}}
\caption{}\l{f19}
\medskip
\end{figure}

In case when $q_1(\Pi) >  q_2(\Pi)$ starting from $\Omega$ first construct a sunflower $\Omega_1$    
such that $${\Omega_1}_1=\{b_{1,l+1}(\Pi), b_{1,l+2}(\Pi), ... b_{1,q_1(\Pi)}(\Pi)\},$$
where $l=q_1(\Pi)-q_2(\Pi),$ 
${\Omega_1}_2=\Pi_2$, and $i(\Pi) = s$ (see again Fig. \ref{f19}).
It is easy to see that for 
the number $\nu$ 
of 1-vertices 
of valency 1 of $\Omega_1$ 
the equality 
$$\nu=\sum_{j=1}^{q_2(\Pi)} (b_{2,j}(\Pi)-2)$$ 
holds (this formula turns out to be true for any choice of the color for the vertex of valency 3 on Fig. \ref{f18}). 
Since by \eqref{eche} \be \nu \geq q_1(\Pi)-q_2(\Pi)\ee
and all vertices of valency 1 of $\Omega_1$ are adjacent to the exterior 
face of $\Omega_1$ by construction, it follows that after gluing a number of edges to the 1-vertices of valency 1 of $\Omega_1$ we obtain 
a sunflower $\Omega_2$ for which 
${\Omega_2}_i=\Pi_i,$ $1\leq i \leq r,$ and $i(\Omega_2)=s$
(see Fig. \ref{f20}).

\begin{figure}[h]
\renewcommand{\captionlabeldelim}{.}
\medskip
\epsfxsize=3truecm
\centerline{\epsffile{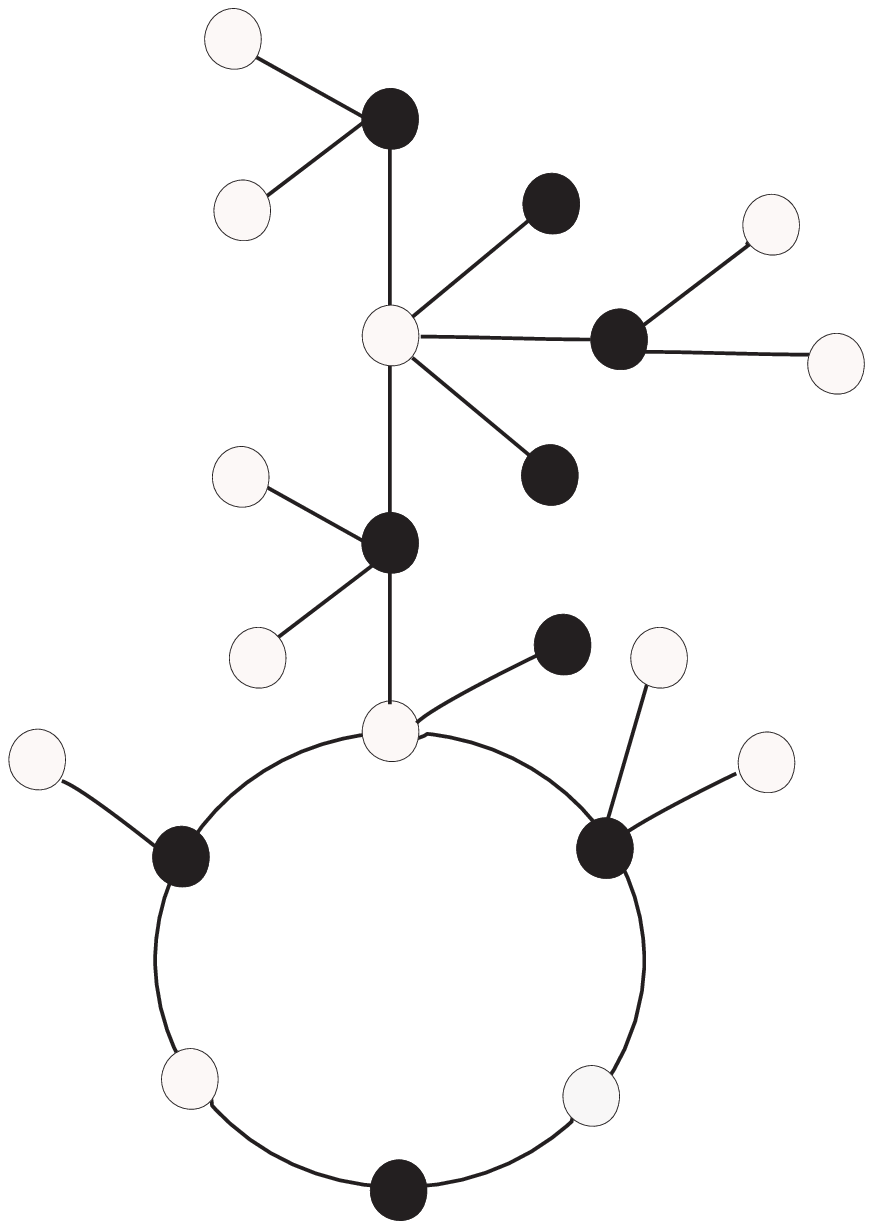}}
\caption{}\l{f20}
\medskip
\end{figure}

For $s=q_2(\Pi)$ the proof of the lemma is similar. The only difference is that we start from the chain $\Omega$ all the vertices of which have valency 2 and
$q_1(\Omega)=q_2(\Omega)=q_2(\Pi),$ $i(\Omega)=s$ (see Fig. \ref{f21}, a, b, c). 

\begin{figure}[h]
\renewcommand{\captionlabeldelim}{.}
\medskip
\epsfxsize=8truecm
\centerline{\epsffile{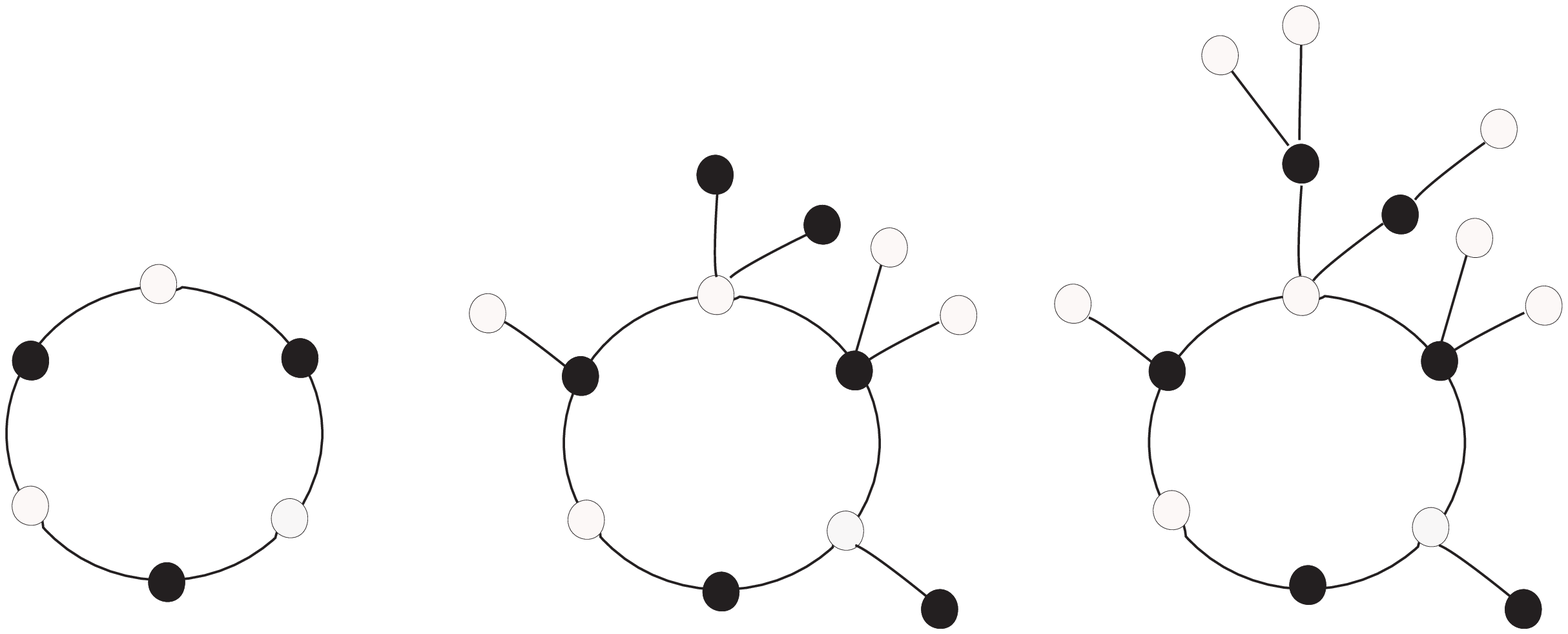}}
\caption{}\l{f21}
\medskip
\end{figure}

\bl \l{pp} Any  passport $\Pi$ with $r(\Pi)=2$ for which  
$q_1(\Pi)=q_2(\Pi)$ 
is realizable whenever $\Pi_1\neq\{2,2,\, ...\, ,2\}$ or $\Pi_2\neq \{2,2,\, ...\, ,2\}.$
\el 

\pr 
If $s\leq q_2(\Pi)$ then the proposition follows from lemma \ref{chai}. 
To prove it for  
$q_2(\Pi) < s \leq n/2$ 
we begin from the sunflower $\Omega$ for which $\Omega_i=\Pi_i$, $1\leq i \leq 2,$ and $i(\Omega)=q_2(\Pi)$
shown on Fig. \ref{f21}, b and then start
shifting its branches from outside to inside. 
Clearly, in this way we can obtain the sunflower $\Omega_1$ with ${\Omega_1}_i=\Pi_i$, $1\leq i \leq 2,$ and $i(\Omega_1)=s$ for any $s$ such that
$q_2(\Pi) < s \leq \mu,$ where $$\mu=q_2(\Pi)+\sum_{j=1}^{q_2(\Pi)} (b_{2,j}(\Pi)-2).$$ Since $\mu$ coincides with the value given by formula \eqref{fo1} for $r=2$ now the lemma follows from formulas \eqref{ffo2}, 
\eqref{fo2}.

\bl \l{aq} Any  passport $\Pi$ with $r(\Pi)=2$ for which $q_1(\Pi)>q_2(\Pi)$ and 
$ q_2(\Pi) < s \leq n/2$ 
is realizable whenever $\Pi_1\neq \{2,2, ... ,2\}$ 
and $q_2(\Pi)>1$.
\el 
 
\pr The proof of the lemma is similar to the one of lemma \ref{ss}. Starting from the sunflower $\Omega$ 
for which $\Omega_i=\Pi_i$, $1\leq i \leq 2,$ and $i(\Omega)=q_2(\Pi)$ shown on Fig. \ref{f21}, c and 
shifting its branches from outside to inside 
we can obtain 
a constellation $\Omega_1$ for which ${\Omega_1}_i=\Pi_i$, $1\leq i \leq 2,$ and $i(\Omega_1)=s,$
where $s$ is any number which can be represented as a sum
\be \l{we} s=q_2(\Pi)+y+x_1b_{1,1}(\Pi)+x_2b_{1,2}(\Pi)+ ... + x_lb_{1,l}(\Pi)\ee for some $x_i,y$ 
with $x_i \in \{0,1\},$ $1\leq i \leq l,$ and $y \in \{0,1,2, ...\, ,t\},$ where 
$l=q_1(\Pi)-q_2(\Pi)$ and \be \l{wer} t= \sum_{j=l+1}^{q_1(\Pi)}(b_{1,j}(\Pi)-2) +\sum_{j=1}^{q_2(\Pi)}(b_{2,j}(\Pi)-2)-l.\ee Formulas \eqref{we}, \eqref{wer} are particular cases for $r=2$ of formulas \eqref{sss}, \eqref{ff}, in particular 
inequality \eqref{fo4} holds for $s_{\rm max}$.

As in lemma \ref{ss} to finish the proof it is enough to establish formula \eqref{los} and for this purpose
it is enough to prove formulas \eqref{fo5}, \eqref{loss}. Formula \eqref{fo5}
now follows directly from the condition $q_2(\Pi)>1$ 
while formula \eqref{loss} follows 
as in lemma \ref{ss} from the condition $\Pi_1\neq \{2,2, ... ,2\}$.
 
\bl \l{iop} Any passport $\Pi$ with $r(\Pi)=2$ for which 
$q_2(\Pi)=1$ 
is realizable whenever $\Pi$ is distinct from  
$\Pi_1=\{l,l,\, ...\, ,l\},$ $\Pi_2=\{1,1,\, ...\, ,1, d\},$ 
$\Pi_3$=\linebreak$\{s,n-s\},$ 
where $l \geq 2,$ $d\geq 3,$ and $s\equiv 0\ \mod l.$     
\el 

\pr It is easy to see that if $\Pi$ 
is realizable then the corresponding constellation
$\Sigma$ has the form shown on Fig. \ref{f22} (1-vertices are colored by the black color).
\begin{figure}[t]
\renewcommand{\captionlabeldelim}{.}
\medskip
\epsfxsize=2truecm
\centerline{\epsffile{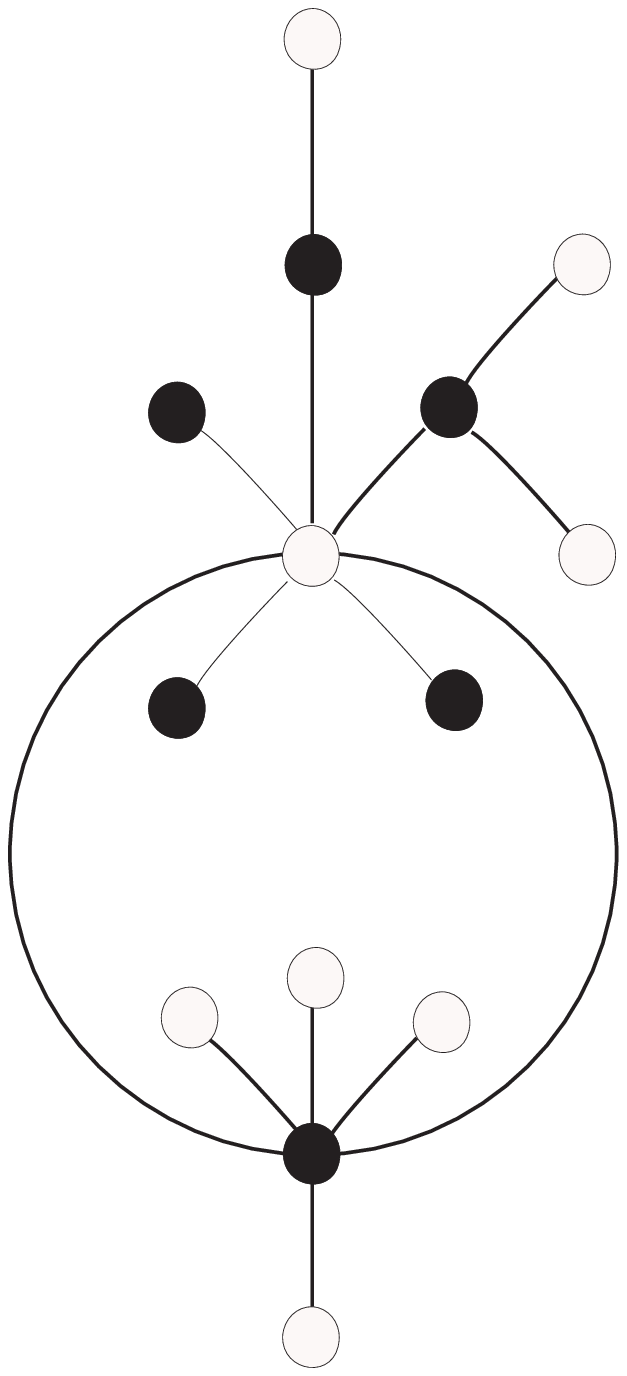}}
\caption{}\l{f22}
\medskip
\end{figure}
Furthermore, we can assume that $b_{2,1}\geq 3$ since otherwise $q_1(\Sigma)=q_2(\Sigma)$ and $\Pi$ is realizable by lemma \ref{pp}.

Placing a 1-vertex of the maximal valency on the cycle and acting as in the proof of lemma \ref{aq} we can obtain a constellation $\Omega$ with 
$\Omega_i=\Pi_i$, $1\leq i \leq 2,$ and $i(\Omega)=s$ for any $s$ which can be represented 
in the form

$$s=1+y+x_1b_{1,1}(\Pi)+x_2b_{1,2}(\Pi)+ ... + x_{q_1-1}b_{1,q_1-1}(\Pi)$$ for some 
$x_i \in \{0,1\},$ $1\leq i \leq q_1(\Pi)-1,$ and $y \in \{0,1,2, ...\, ,t\},$
where $$t=b_{1,q_1}(\Pi)-2+(b_{2,1}(\Pi)-2)-(q_1(\Pi)-1)$$ 
(these formulas are particular cases for $q_2(\Pi)=1$ of formulas \eqref{we}, \eqref{wer} in particular 
inequality \eqref{fo4} holds for $s_{\rm max}$).

Observe that by formula \eqref{eche}
$$t=b_{1,q_1}(\Pi)-2+e_1(\Pi).$$ Therefore, if $e_1(\Pi)>0$ then 
formula \eqref{los} holds
and as above lemma \ref{cal} implies that $\Pi$ is realizable. 

Similarly, if $b_{1,1}(\Pi)< b_{1,q_1}(\Pi)$ then
$\Pi$ is also realizable since in this case all conditions \eqref{us} hold. 
Indeed, for $k=1$ we have  $$t\geq b_{1,q_1}(\Pi)-2\geq b_{1,1}(\Pi)-1$$ and for $k>1$ we have: 
$$t+\sum_{i=1}^{k-1}b_{1,i}(\Pi)\geq t+2 \geq b_{1,q_1}(\Pi) \geq b_{1,k}(\Pi)>b_{1,k}(\Pi)-1.$$

It follows that $\Pi$ may not be realizable only if $e_1=0$ and $b_{1,1}= b_{1,q_1}$ that is
if $\Pi_1=\{l,l,\, ...\, ,l\},$ $\Pi_2=\{1,1,\, ...\, ,1, d\}$. 
Now it is easy to establish by a direct calculation that such $\Pi$ is realizable if and only if
$s\not\equiv 0\ \mod l.$

\bl \l{zxc} Let $\Pi$ be a passport with $r(\Pi)=2$ for which $\Pi_1= \{2,2, ... ,2\},$ 
$\Pi_2\neq \{2,2, ... ,2\}$ and $q_2(\Pi)>1.$ Suppose that \be \l{neus} \sum_{i=2}^{q_2(\Pi)} (b_{2,i}(\Pi)-2)< b_{2,1}.\ee
Then either 
\vskip 0.2cm
\noindent 1)
$\Pi_2=\{1,1,\, ... \, , 1,d,d\},$ where $d\geq 3$, or 
\vskip 0.2cm
\noindent 2)
$\Pi_2=\{1,1,\, ... \, , 1,d-1,d\},$ where $d\geq 3$, or  
\vskip 0.2cm
\noindent 3)
$\Pi_2=\{1,1,1,3,3,3\},$ or
\vskip 0.2cm
\noindent 4)
$\Pi_2=\{1,2,2,\, ... \, ,2 
,3\}.$
\el  

\pr If $q_2(\Pi)=2$ then
$$\sum_{i=2}^{q_2(\Pi)} (b_{2,i}(\Pi)-2)= b_{2,2}(\Pi)-2$$
and \eqref{neus} holds only if  
$$\Pi_2=\{1,1,\, ... \, , 1,d,d\} \ \ \ \ {\rm or} \ \ \ \ 
\Pi_2=\{1,1,\, ... \, , 1,d-1,d\},$$ where $d=b_{2,q_2}(\Pi)\geq 3$ in view of lemma \ref{ll}. So, in the following we will 
assume that $q_2(\Pi)\geq 3$.

If $b_{2,1}(\Pi)\geq 3$ then
$$\sum_{i=2}^{q_2(\Pi)} (b_{2,i}(\Pi)-2)\geq b_{2,2}(\Pi)+b_{2,3}(\Pi)-4\geq  
2b_{2,1}(\Pi)-4\geq
b_{2,1}(\Pi)-1$$ with the equality only if $q_2(\Pi)=3$ and $b_{2,3}=b_{2,2}=b_{2,1}=3.$ 
Therefore, in this case 
condition \eqref{neus} holds 
only if $\Pi_2=\{1,1,\, ... \, , 1,3,3,3\}.$ Denoting the number of appearances of the unit in $\Pi_2$ 
by $l_1$ we obtain that $l_1+9=n$ and $l_1+3=n/2.$ It follows that $l_1=3.$

Suppose now that $b_{2,1}(\Pi)=2.$ In view of lemma \ref{ll} we have $b_{2,q_2}(\Pi)>2.$ If $b_{2,q_2}(\Pi)>3$ then 
$$\sum_{i=2}^{q_2(\Pi)} (b_{2,i}(\Pi)-2)\geq 2(b_{2,q_2}(\Pi)-2)\geq 4 >b_{2,1}(\Pi).$$
On the other hand, if $b_{2,q_2}(\Pi)=3$ and $b_{2,q_2(\Pi)-1}=3$ then 
$$\sum_{i=2}^{q_2(\Pi)} (b_{2,i}(\Pi)-2)\geq  2(b_{2,q_2}(\Pi)-2)\geq 2=b_{2,1}(\Pi).$$ Hence,
\eqref{neus} holds only if $b_{2,q_2(\Pi)}=3,$ $b_{2,q_2(\Pi)-1}=2$ or equivalently
if $\Pi_2=\{1,1,\, ... \, ,1, 2,2,\, ... \, ,2 
,3\}.$ Denoting the number of appearances of the number $i,$ $1\leq i \leq 2,$ in $\Pi_2$ 
by $l_i$ we obtain that $l_1+2l_2+3=n$ and $l_1+l_2+1=n/2.$ It follows that $l_1=1.$ 
 
\bl \l{poi} Any passport $\Pi$ with $r(\Pi)=2$ for which $\Pi_1= \{2,2, ... ,2\},$ 
$\Pi_2\neq \{2,2, ... ,2\},$ and $q_2(\Pi)>1$ is realizable whenever 
$\Pi$ is distinct from the passports listed below:
\vskip 0.2cm
\noindent 1)
$\Pi_1=\{2,2, \, ... \, , 2\},$ $\Pi_2=\{1,1,\, ... \, , 1,d,d\},$ 
$\Pi_3=\{2d-3,n-2d+3\},$ 
\vskip 0.2cm
\noindent 2)
$\Pi_1=\{2,2, \, ... \, , 2\},$ $\Pi_2=\{1,1,\, ... \, , 1,d,d\},$ 
$\Pi_3=\{2d-1,n-2d+1\},$ 
\vskip 0.2cm
\noindent 3)
$\Pi_1=\{2,2, \, ... \, , 2\},$ $\Pi_2=\{1,1,\, ... \, , 1,d-1,d\},$ 
$\Pi_3=\{2d-3,n-2d+3\},$  
\vskip 0.2cm
\noindent 4)
$\Pi_1=\{2,2,2,2,2,2\},$ $\Pi_2=\{1,1,1,3,3,3\},$ 
$\Pi_3=\{6,6\},$   
\vskip 0.2cm
\noindent 5)
$\Pi_1=\{2,2, \, ... \, , 2\},$ $\Pi_2=\{1,2,2,\, ... \, ,2 
,3\},$  $\Pi_3=\{n/2,n/2\},$ 
\vskip 0.2cm
\noindent
where $d\geq 3.$  
\el  
 
\pr In view of lemma \ref{chai} if $s(\Pi)\leq q_2(\Pi)$ then $\Pi$ is realizable 
so we only must consider the case when $q_2(\Pi)< s(\Pi) \leq  n/2.$ 

First observe that if $s\equiv q_2(\Pi) \ \mod 2$ then $\Pi$ is realizable. 
Indeed, starting from a constellation $\Gamma$ for which $\Gamma_1=\Pi_1,$ 
$\Gamma_2=\Pi_2,$ and $i(\Gamma)=q_2(\Pi)$ constructed in lemma \ref{chai}
(see  Fig. \ref{f24},a, where the condition $\Pi_1=\{2,2, \, ... \, , 2\}$ is reflected)  
\begin{figure}[h]
\renewcommand{\captionlabeldelim}{.}
\medskip
\epsfxsize=9truecm
\centerline{\epsffile{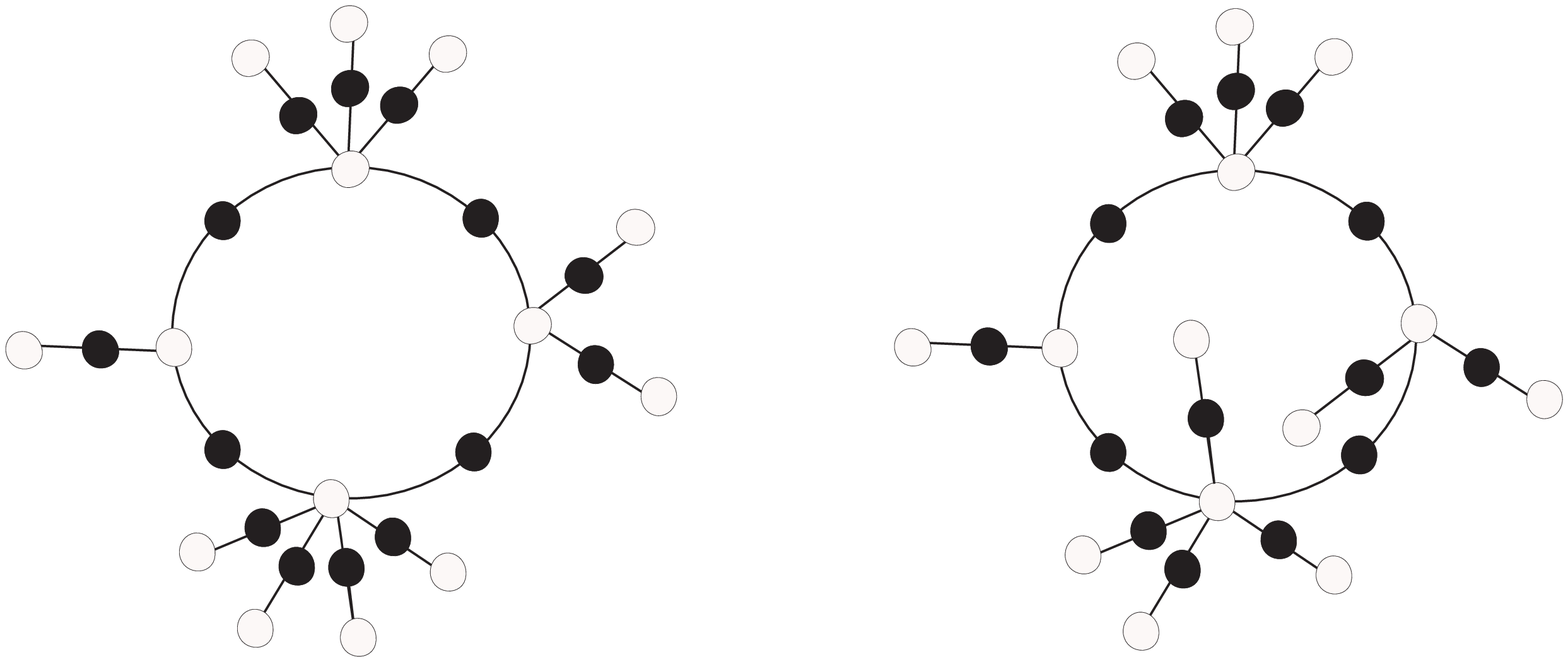}}
\caption{}\l{f24}
\medskip
\end{figure}
and shifting the branches of $\Gamma$
from outside to inside (see Fig. \ref{f24},b) one can obtain
a constellation $\Sigma$ for which $\Sigma_1=\Pi_1,$ 
$\Sigma_2=\Pi_2,$
and $i(\Sigma)=s$ for any $s\equiv q_2(\Pi)\ \mod 2$ such that $$s\leq q_2(\Pi) +2\sum_{j=1}^{q_2(\Pi)} (b_{2,j}(\Pi)-2).$$
Since in view of \eqref{eche} $$s_{\max}=q_2(\Pi)+2(e_1(\Pi) +q_1(\Pi)-q_2(\Pi))=2(e_1(\Pi) +q_1(\Pi))-q_2(\Pi)=
n-q_2(\Pi)$$ and $n-q_2(\Pi)\geq n/2$
this implies the statement.
 
Consider now the case when $s\equiv 1+q_2(\Pi) \ \mod 2$. 
Modify the constellation shown on Fig. \ref{f24},a so that to obtain   
a constellation $\tilde \Gamma$ for which all 2-vertices of valency $>1$ except one are on the cycle (see  Fig. \ref{f25}, a) and the valency of the exceptional vertex is $b_{2,1}$
(recall that $q_2(\Pi)\geq 2$ and that $b_{2,q_2}(\Pi)>2$ by lemma \ref{ll}). Clearly, we have    
$\tilde \Gamma_1=\Pi_1,$ $\tilde \Gamma_2=\Pi_2,$
and $i(\tilde \Gamma)=q_2(\Pi)-1$.  
\begin{figure}[h]
\renewcommand{\captionlabeldelim}{.}
\medskip
\epsfxsize=8truecm
\centerline{\epsffile{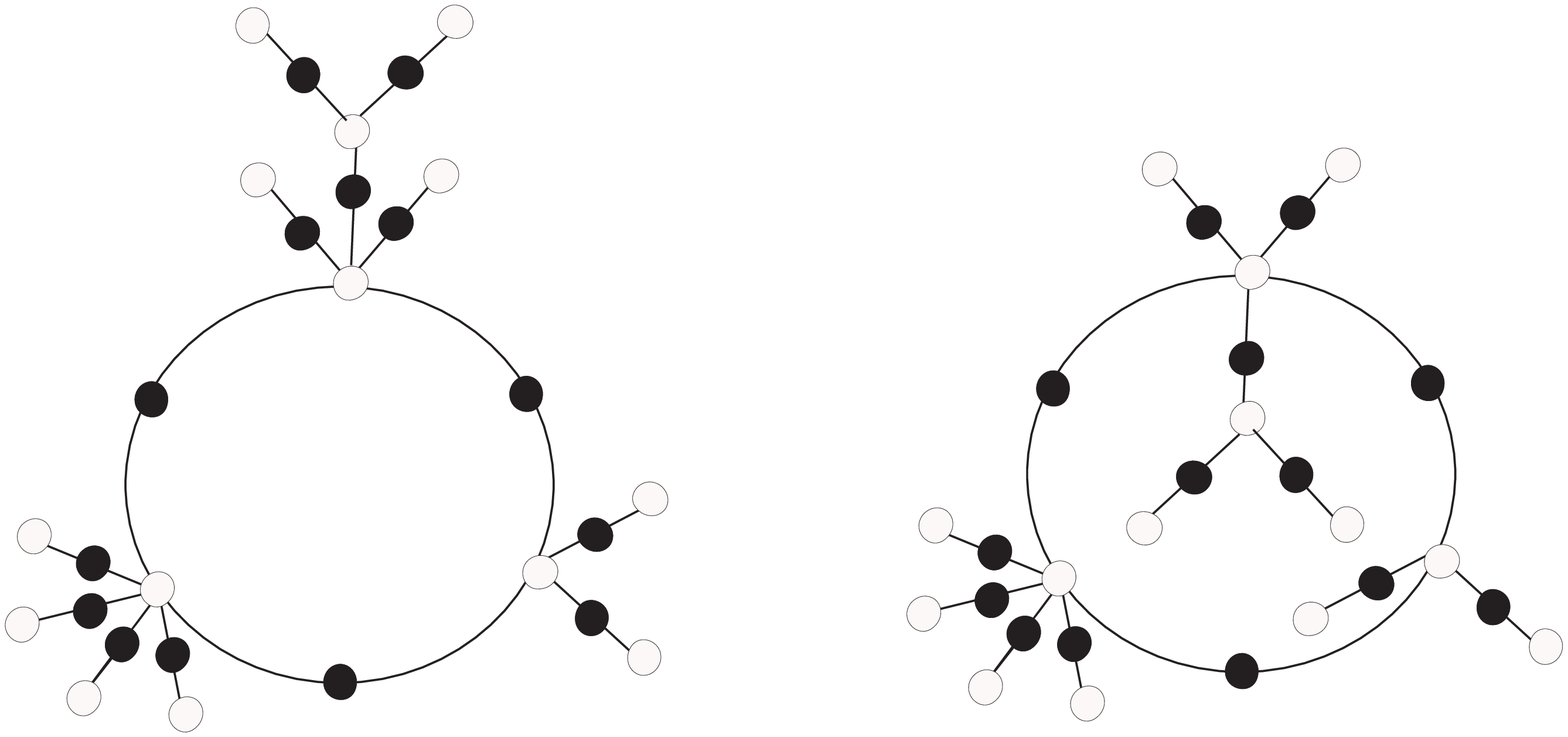}}
\caption{}\l{f25}
\medskip
\end{figure}
Shifting now the branches of $\tilde \Gamma$
from outside to inside (see Fig. \ref{f25}, b) one can obtain
a constellation $\Sigma$ for which $\Sigma_1=\Pi_1,$ 
$\Sigma_2=\Pi_2,$
and $i(\Sigma)=s$ for any $s\equiv 1+q_2(\Pi)\ \mod 2$ which can be represented as 
$$s=q_2(\Pi)-1 +2y+2b_{2,1}x$$ with $x\in \{0,1\}$ and $y \in \{0,1,\, ...\, ,t\},$ where   
$$t= \sum_{i=2}^{q_2(\Pi)} (b_{2,i}(\Pi)-2)-1.$$

Furthermore, in view of lemma \ref{cal} if \be \l{uss}
\sum_{i=2}^{q_2(\Pi)} (b_{2,i}(\Pi)-2)\geq b_{2,1} \ee
then we obtain in this way any 
$s$ such that $$s\equiv 1+q_2(\Pi)\ \mod 2, \ \ \ \ \ \ \  \ \ q_2(\Pi) < s \leq s_{\max}.$$ 
Since in view of \eqref{eche}
$$s_{\max}=q_2(\Pi)-1+2(e_1(\Pi)+q_1(\Pi)-q_2(\Pi)-(b_{2,1}(\Pi)-2))-2+
2b_{2,1}(\Pi)=$$
$$=-q_2(\Pi)+1+2(e_1(\Pi)+q_1(\Pi))=
n-q_2(\Pi)+1\geq  n/2$$ it follows that in order
to prove the lemma we only must investigate when the passports listed in lemma \ref{zxc} are realizable.

First of all observe that if for some constellation $\Gamma$ we have:
$\Gamma_1=\{2,2,2,2,2,2\},$ $\Gamma_2=\{1,1,1,3,3,3\}$ and 
\be \l{pppp} i(\Gamma)\equiv 1+q_2(\Pi)\equiv 
0 \ \mod 2, \ \ \ \ \ \ q_2(\Pi)=3< i(\Gamma)\leq n/2=6\ee
then the first of conditions \eqref{pppp} together with the condition 
$\Gamma_1=\{2,2,2,2,2,2\}$ imply that 
the cycle of $\Gamma$ can contain only an even number of 2-vertices. Therefore this number equals 2 and it is easy to see
that $\Gamma$ necessarily has the form
shown on Fig. \ref{f26}.
It follows that a passport $\Pi$ for which $\Pi_1=\{2,2,2,2,2,2\},$ $\Pi_2=\{1,1,1,3,3,3\}$ 
is realizable whenever $\Pi_3$ is distinct from $\{6,6\}.$
\begin{figure}[h]
\renewcommand{\captionlabeldelim}{.}
\medskip
\epsfxsize=2.5truecm
\centerline{\epsffile{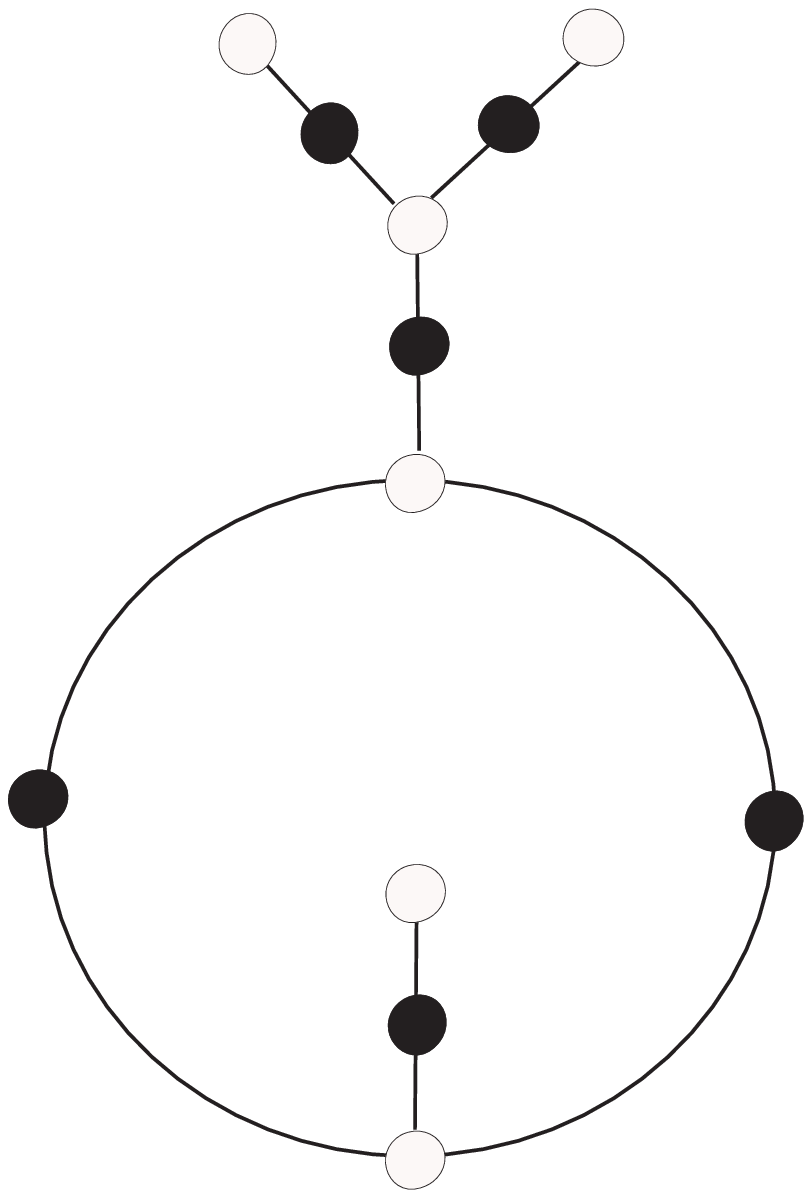}}
\caption{}\l{f26}
\medskip
\end{figure}

Furthermore, if for some constellation $\Gamma$ we have
$\Gamma_1=\{2,2, \, ... \, , 2\},$ $\Gamma_2
=\{1,2,2,\, ... \, ,2 
,3\}$
then it is easy to see that $\Gamma$ has the form shown on Fig. \ref{f28}.
Moreover, since for such $\Gamma$ the equality $q_2(\Gamma)=n/2-1$ holds, the 
condition $q_2(\Pi)< i(\Gamma)\leq n/2$ turns out to be equivalent to the condition
$i(\Gamma)= n/2.$ Clearly, this condition can not be realized for such $\Gamma$ and therefore
a passport $\Pi$ for which  
$\Pi_1=\{2,2, \, ... \, , 2\},$ $\Pi_2=\{1,2,2,\, ... \, ,2 
,3\}$ is realizable whenever  
$\Pi_3\neq\{n/2,n/2\}.$  
\begin{figure}[h]
\renewcommand{\captionlabeldelim}{.}
\medskip
\epsfxsize=2.5truecm
\centerline{\epsffile{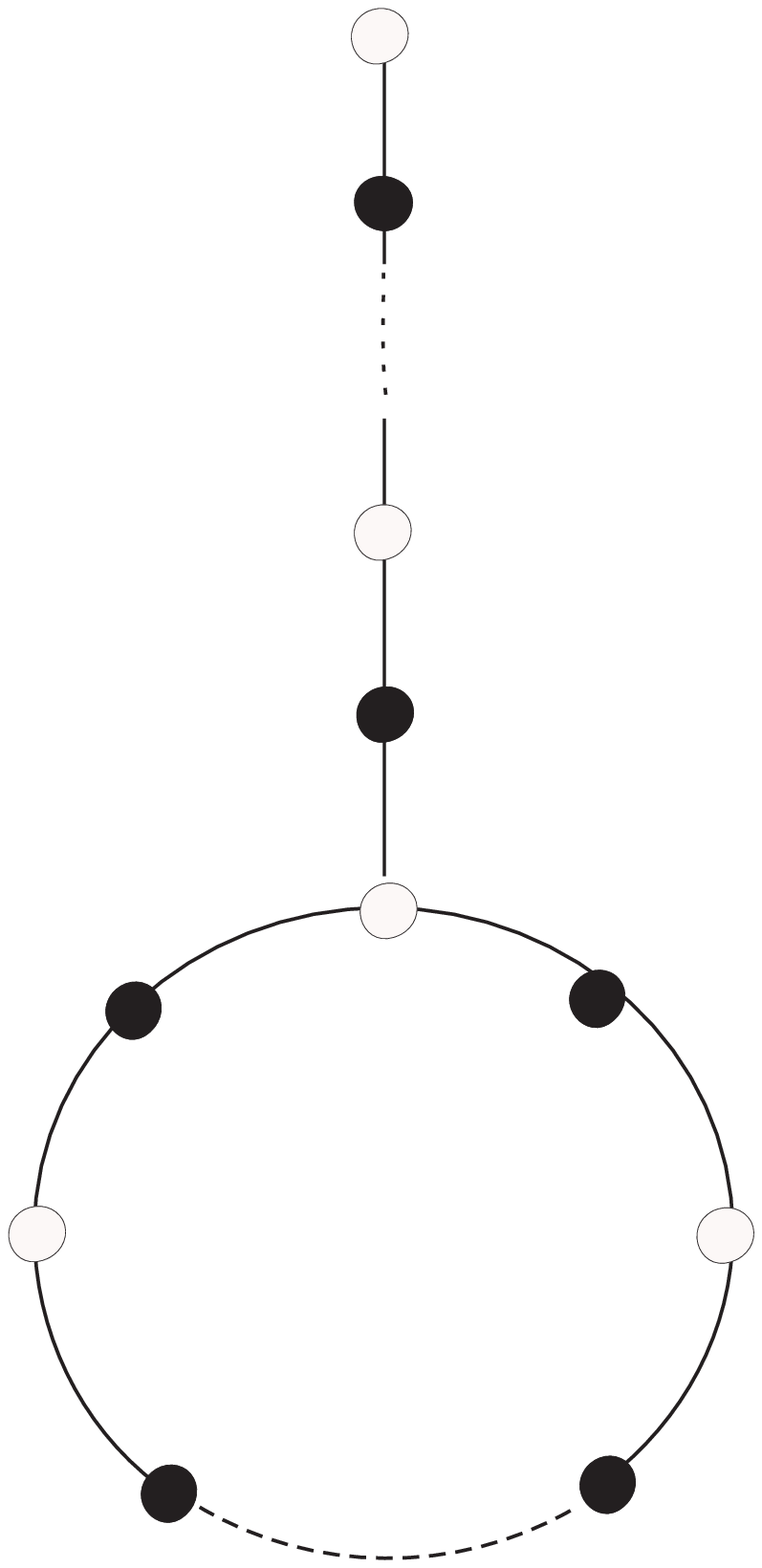}}
\caption{}\l{f28}
\medskip
\end{figure}

Finally if for a constellation $\Gamma$ we have    
$$\Gamma_1=\{2,2, \, ... \, , 2\}, \ \ \ \ \Gamma_2=\{1,1,\, ... \, , 1,d-1,d\}, \ \ \ \ i(\Gamma)\equiv 1+q_2(\Pi)\equiv 
1 \ \mod 2$$ then the cycle of $\Gamma$ contains only one 2-vertex
which is of valency $d$ or of valency $d-1$ and therefore
$\Gamma$ necessarily has the form
shown on Fig. \ref{f27}, a or b.
\begin{figure}[h]
\renewcommand{\captionlabeldelim}{.}
\medskip
\epsfxsize=6truecm
\centerline{\epsffile{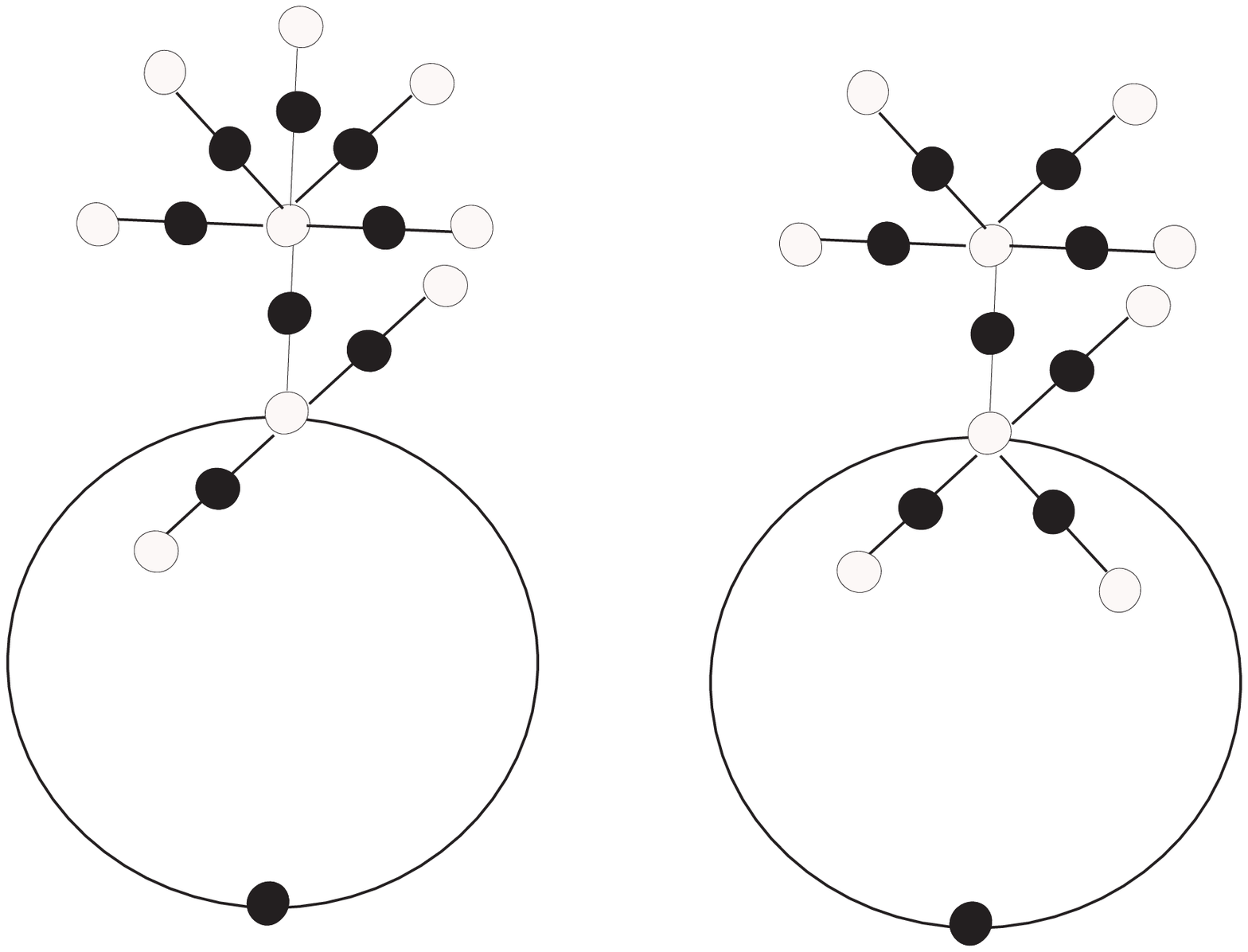}}
\caption{}\l{f27}
\medskip
\end{figure}
It follows easily that a passport $\Pi$ for which $\Pi_1=\{2,2, \, ... \, , 2\},$ 
$\Pi_2=\{1,1,\, ... \, , 1,d-1,d\}$ is realizable whenever $\Pi_3$ distinct from $\Pi_3=\{2d-3,n-2d+3\}.$ 

In the same way one can show that a passport $\Pi$ for which $\Pi_1=\{2,2, \, ... \, , 2\},$ $\Pi_2=\{1,1,\, ... \, , 1,d,d\}$ is realizable whenever 
$\Pi_3\neq\{2d-3,n-2d+3\},$
$\Pi_3\neq\{2d-1,n-2d+1\}.$ 

\bt A passport with $r(\Pi)=2$ is realizable whenever $\Pi$ is distinct from the passports
listed in the main theorem.
\et  

\pr Indeed, if a passport $\Pi$ with $\Pi_1=\{2,2,2, ... 2\},$ $\Pi_2=\{2,2,2, ... 2\}$ is realizable then 
the bicolored graph $\Gamma$ corresponding to $\Pi$ should have the form shown on Fig. \ref{c5} and therefore $s= n/2.$
\begin{figure}[h]
\renewcommand{\captionlabeldelim}{.}
\medskip
\epsfxsize=2.5truecm
\centerline{\epsffile{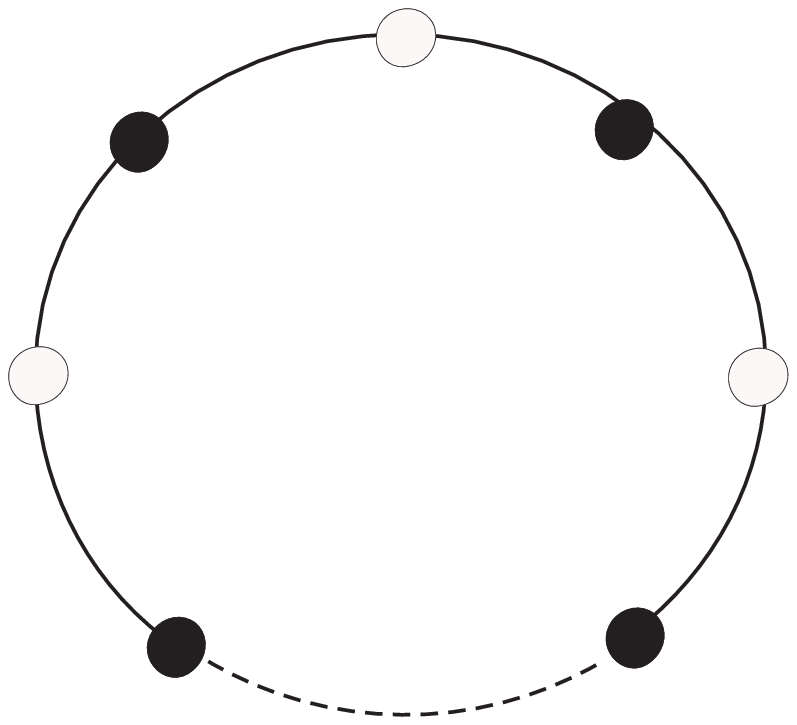}}
\caption{}\l{c5}
\medskip
\end{figure}  
So, we can assume that
either $\Pi_1\neq \{2,2, ... ,2\}$ or $\Pi_2\neq \{2,2, ... ,2\}$.
Furthermore, in view of lemmas \ref{chai}, \ref{pp}, \ref{aq} such a passport
may not be realizable only if $\Pi_1= \{2,2, ... ,2\}$ or $q_2(\Pi)=1$.

If $q_2(\Pi)=1$ then by lemma \ref{iop} the passport $\Pi$ is realizable whenever 
it is distinct from the passport 1).
On the other hand, 
if $q_2(\Pi)>1$ but $\Pi_1= \{2,2, ... ,2\}$ then by lemma \ref{poi} the passport $\Pi$ is realizable whenever 
it is distinct from the passports 3)-7). 

\end{document}